
\ifx\shlhetal\undefinedcontrolsequence\let\shlhetal\relax\fi

\input amstex
\expandafter\ifx\csname mathdefs.tex\endcsname\relax
  \expandafter\gdef\csname mathdefs.tex\endcsname{}
\else \message{Hey!  Apparently you were trying to
  \string\input{mathdefs.tex} twice.   This does not make sense.} 
\errmessage{Please edit your file (probably \jobname.tex) and remove
any duplicate ``\string\input'' lines}\endinput\fi




\catcode`\X=12\catcode`\@=11

\def\n@wcount{\alloc@0\count\countdef\insc@unt}
\def\n@wwrite{\alloc@7\write\chardef\sixt@@n}
\def\n@wread{\alloc@6\read\chardef\sixt@@n}
\def\r@s@t{\relax}\def\v@idline{\par}\def\@mputate#1/{#1}
\def\l@c@l#1X{\firstpart.#1}\def\gl@b@l#1X{#1}\def\t@d@l#1X{{}}

\def\crossrefs#1{\ifx\all#1\let\tr@ce=\all\else\def\tr@ce{#1,}\fi
   \n@wwrite\cit@tionsout\openout\cit@tionsout=\jobname.cit 
   \write\cit@tionsout{\tr@ce}\expandafter\setfl@gs\tr@ce,}
\def\setfl@gs#1,{\def\@{#1}\ifx\@\empty\let\next=\relax
   \else\let\next=\setfl@gs\expandafter\xdef
   \csname#1tr@cetrue\endcsname{}\fi\next}
\def\m@ketag#1#2{\expandafter\n@wcount\csname#2tagno\endcsname
     \csname#2tagno\endcsname=0\let\tail=\all\xdef\all{\tail#2,}
   \ifx#1\l@c@l\let\tail=\r@s@t\xdef\r@s@t{\csname#2tagno\endcsname=0\tail}\fi
   \expandafter\gdef\csname#2cite\endcsname##1{\expandafter
     \ifx\csname#2tag##1\endcsname\relax?\else\csname#2tag##1\endcsname\fi
     \expandafter\ifx\csname#2tr@cetrue\endcsname\relax\else
     \write\cit@tionsout{#2tag ##1 cited on page \folio.}\fi}
   \expandafter\gdef\csname#2page\endcsname##1{\expandafter
     \ifx\csname#2page##1\endcsname\relax?\else\csname#2page##1\endcsname\fi
     \expandafter\ifx\csname#2tr@cetrue\endcsname\relax\else
     \write\cit@tionsout{#2tag ##1 cited on page \folio.}\fi}
   \expandafter\gdef\csname#2tag\endcsname##1{\expandafter
      \ifx\csname#2check##1\endcsname\relax
      \expandafter\xdef\csname#2check##1\endcsname{}%
      \else\immediate\write16{Warning: #2tag ##1 used more than once.}\fi
      \multit@g{#1}{#2}##1/X%
      \write\t@gsout{#2tag ##1 assigned number \csname#2tag##1\endcsname\space
      on page \number\count0.}%
   \csname#2tag##1\endcsname}}

\def\multit@g#1#2#3/#4X{\def\t@mp{#4}\ifx\t@mp\empty%
      \global\advance\csname#2tagno\endcsname by 1 
      \expandafter\xdef\csname#2tag#3\endcsname
      {#1\number\csname#2tagno\endcsnameX}%
   \else\expandafter\ifx\csname#2last#3\endcsname\relax
      \expandafter\n@wcount\csname#2last#3\endcsname
      \global\advance\csname#2tagno\endcsname by 1 
      \expandafter\xdef\csname#2tag#3\endcsname
      {#1\number\csname#2tagno\endcsnameX}
      \write\t@gsout{#2tag #3 assigned number \csname#2tag#3\endcsname\space
      on page \number\count0.}\fi
   \global\advance\csname#2last#3\endcsname by 1
   \def\t@mp{\expandafter\xdef\csname#2tag#3/}%
   \expandafter\t@mp\@mputate#4\endcsname
   {\csname#2tag#3\endcsname\lastpart{\csname#2last#3\endcsname}}\fi}
\def\t@gs#1{\def\all{}\m@ketag#1e\m@ketag#1s\m@ketag\t@d@l p
\let\realscite\scite
\let\realstag\stag
   \m@ketag\gl@b@l r \n@wread\t@gsin
   \openin\t@gsin=\jobname.tgs \re@der \closein\t@gsin
   \n@wwrite\t@gsout\openout\t@gsout=\jobname.tgs }
\outer\def\localtags{\t@gs\l@c@l}
\outer\def\globaltags{\t@gs\gl@b@l}
\outer\def\newlocaltag#1{\m@ketag\l@c@l{#1}}
\outer\def\newglobaltag#1{\m@ketag\gl@b@l{#1}}

\newif\ifpr@ 
\def\m@kecs #1tag #2 assigned number #3 on page #4.%
   {\expandafter\gdef\csname#1tag#2\endcsname{#3}
   \expandafter\gdef\csname#1page#2\endcsname{#4}
   \ifpr@\expandafter\xdef\csname#1check#2\endcsname{}\fi}
\def\re@der{\ifeof\t@gsin\let\next=\relax\else
   \read\t@gsin to\t@gline\ifx\t@gline\v@idline\else
   \expandafter\m@kecs \t@gline\fi\let \next=\re@der\fi\next}
\def\pretags#1{\pr@true\pret@gs#1,,}
\def\pret@gs#1,{\def\@{#1}\ifx\@\empty\let\n@xtfile=\relax
   \else\let\n@xtfile=\pret@gs \openin\t@gsin=#1.tgs \message{#1} \re@der 
   \closein\t@gsin\fi \n@xtfile}

\newcount\sectno\sectno=0\newcount\subsectno\subsectno=0
\newif\ifultr@local \def\ultralocal{\ultr@localtrue}
\def\firstpart{\number\sectno}
\def\lastpart#1{\ifcase#1 \or a\or b\or c\or d\or e\or f\or g\or h\or 
   i\or k\or l\or m\or n\or o\or p\or q\or r\or s\or t\or u\or v\or w\or 
   x\or y\or z \fi}

\def\resetall{\global\advance\sectno by 1\subsectno=0
   \gdef\firstpart{\number\sectno}\r@s@t}
\def\resetsub{\global\advance\subsectno by 1
   \gdef\firstpart{\number\sectno.\number\subsectno}\r@s@t}
\def\newsection#1\par{\resetall\vskip0pt plus.3\vsize\penalty-250
   \vskip0pt plus-.3\vsize\bigskip\bigskip
   \message{#1}\leftline{\bf#1}\nobreak\bigskip}
\def\subsection#1\par{\ifultr@local\resetsub\fi
   \vskip0pt plus.2\vsize\penalty-250\vskip0pt plus-.2\vsize
   \bigskip\smallskip\message{#1}\leftline{\bf#1}\nobreak\medskip}


\newdimen\marginshift

\newdimen\margindelta
\newdimen\marginmax
\newdimen\marginmin

\def\margininit{       
\marginmax=3 true cm                  
				      
\margindelta=0.1 true cm              
\marginmin=0.1true cm                 
\marginshift=\marginmin
}    

\def\t@gsjj#1,{\def\@{#1}\ifx\@\empty\let\next=\relax\else\let\next=\t@gsjj
   \def\@@{p}\ifx\@\@@\else
   \expandafter\gdef\csname#1cite\endcsname##1{\citejj{##1}}
   \expandafter\gdef\csname#1page\endcsname##1{?}
   \expandafter\gdef\csname#1tag\endcsname##1{\tagjj{##1}}\fi\fi\next}
\newif\ifshowstuffinmargin
\showstuffinmarginfalse
\def\jjtags{\ifx\shlhetal\relax 
  \else
\ifx\shlhetal\undefinedcontrolseq
\else
\showstuffinmargintrue
\ifx\all\relax\else\expandafter\t@gsjj\all,\fi\fi \fi
}

\def\tagjj#1{\realstag{#1}\oldmginpar{\zeigen{#1}}}
\def\citejj#1{\rechnen{#1}\mginpar{\zeigen{#1}}}     

\def\rechnen#1{\expandafter\ifx\csname stag#1\endcsname\relax ??\else
                           \csname stag#1\endcsname\fi}

\newdimen\theight

\def\marginfont{\sevenrm}

\def\trymarginbox#1{\setbox0=\hbox{\marginfont\hskip\marginshift #1}%
		\global\marginshift\wd0 
		\global\advance\marginshift\margindelta}

\def \oldmginpar#1{%
\ifvmode\setbox0\hbox to \hsize{\hfill\rlap{\marginfont\quad#1}}%
\ht0 0cm
\dp0 0cm
\box0\vskip-\baselineskip
\else 
             \vadjust{\trymarginbox{#1}%
		\ifdim\marginshift>\marginmax \global\marginshift\marginmin
			\trymarginbox{#1}%
                \fi
             \theight=\ht0
             \advance\theight by \dp0    \advance\theight by \lineskip
             \kern -\theight \vbox to \theight{\rightline{\rlap{\box0}}%
\vss}}\fi}

\newdimen\upordown
\global\upordown=8pt
\font\tinyfont=cmtt8 
\def\mginpar#1{\smash{\hbox to 0cm{\kern-10pt\raise7pt\hbox{\tinyfont #1}\hss}}}
\def\mginpar#1{{\hbox to 0cm{\kern-10pt\raise\upordown\hbox{\tinyfont #1}\hss}}\global\upordown-\upordown}


\def\t@gsoff#1,{\def\@{#1}\ifx\@\empty\let\next=\relax\else\let\next=\t@gsoff
   \def\@@{p}\ifx\@\@@\else
   \expandafter\gdef\csname#1cite\endcsname##1{\zeigen{##1}}
   \expandafter\gdef\csname#1page\endcsname##1{?}
   \expandafter\gdef\csname#1tag\endcsname##1{\zeigen{##1}}\fi\fi\next}
\def\verbatimtags{\showstuffinmarginfalse
\ifx\all\relax\else\expandafter\t@gsoff\all,\fi}
\def\zeigen#1{\hbox{$\scriptstyle\langle$}#1\hbox{$\scriptstyle\rangle$}}


\def\(#1){\edef\dot@g{\ifmmode\ifinner(\hbox{\noexpand\etag{#1}})
   \else\noexpand\eqno(\hbox{\noexpand\etag{#1}})\fi
   \else(\noexpand\ecite{#1})\fi}\dot@g}

\newif\ifbr@ck
\def\eat#1{}
\def\[#1]{\br@cktrue[\br@cket#1'X]}
\def\br@cket#1'#2X{\def\temp{#2}\ifx\temp\empty\let\next\eat
   \else\let\next\br@cket\fi
   \ifbr@ck\br@ckfalse\br@ck@t#1,X\else\br@cktrue#1\fi\next#2X}
\def\br@ck@t#1,#2X{\def\temp{#2}\ifx\temp\empty\let\neext\eat
   \else\let\neext\br@ck@t\def\temp{,}\fi
   \def\teemp{#1}\ifx\teemp\empty\else\rcite{#1}\fi\temp\neext#2X}
\def\resetbr@cket{\gdef\[##1]{[\rtag{##1}]}}
\def\references{\resetbr@cket\newsection References\par}

\newtoks\symb@ls\newtoks\s@mb@ls\newtoks\p@gelist\n@wcount\ftn@mber
    \ftn@mber=1\newif\ifftn@mbers\ftn@mbersfalse\newif\ifbyp@ge\byp@gefalse
\def\defm@rk{\ifftn@mbers\n@mberm@rk\else\symb@lm@rk\fi}
\def\n@mberm@rk{\xdef\m@rk{{\the\ftn@mber}}%
    \global\advance\ftn@mber by 1 }
\def\rot@te#1{\let\temp=#1\global#1=\expandafter\r@t@te\the\temp,X}
\def\r@t@te#1,#2X{{#2#1}\xdef\m@rk{{#1}}}
\def\b@@st#1{{$^{#1}$}}\def\str@p#1{#1}
\def\symb@lm@rk{\ifbyp@ge\rot@te\p@gelist\ifnum\expandafter\str@p\m@rk=1 
    \s@mb@ls=\symb@ls\fi\write\f@nsout{\number\count0}\fi \rot@te\s@mb@ls}
\def\byp@ge{\byp@getrue\n@wwrite\f@nsin\openin\f@nsin=\jobname.fns 
    \n@wcount\currentp@ge\currentp@ge=0\p@gelist={0}
    \re@dfns\closein\f@nsin\rot@te\p@gelist
    \n@wread\f@nsout\openout\f@nsout=\jobname.fns }
\def\m@kelist#1X#2{{#1,#2}}
\def\re@dfns{\ifeof\f@nsin\let\next=\relax\else\read\f@nsin to \f@nline
    \ifx\f@nline\v@idline\else\let\t@mplist=\p@gelist
    \ifnum\currentp@ge=\f@nline
    \global\p@gelist=\expandafter\m@kelist\the\t@mplistX0
    \else\currentp@ge=\f@nline
    \global\p@gelist=\expandafter\m@kelist\the\t@mplistX1\fi\fi
    \let\next=\re@dfns\fi\next}
\def\symbols#1{\symb@ls={#1}\s@mb@ls=\symb@ls} 
\def\bigsymbol{\textstyle}
\symbols{\bigsymbol\ast,\dagger,\ddagger,\sharp,\flat,\natural,\star}
\def\ftnumbers{\ftn@mberstrue} \def\ftsymbols{\ftn@mbersfalse}
\def\paginal{\byp@ge} \def\resetftnumbers{\ftn@mber=1}
\def\ftnote#1{\defm@rk\expandafter\expandafter\expandafter\footnote
    \expandafter\b@@st\m@rk{#1}}

\long\def\jump#1\endjump{}
\def\ssum{\mathop{\lower .1em\hbox{$\textstyle\Sigma$}}\nolimits}

\def\qed{\nobreak\kern 1em \vrule height .5em width .5em depth 0em}
\def\newneq{\hbox{\rlap{\hbox to 1\wd9{\hss$=$\hss}}\raise .1em 
   \hbox to 1\wd9{\hss$\scriptscriptstyle/$\hss}}}
\def\subsetne{\setbox9 = \hbox{$\subset$}\mathrel{\hbox{\rlap
   {\lower .4em \newneq}\raise .13em \hbox{$\subset$}}}}
\def\supsetne{\setbox9 = \hbox{$\subset$}\mathrel{\hbox{\rlap
   {\lower .4em \newneq}\raise .13em \hbox{$\supset$}}}}

\def\vbar{\mathchoice{\vrule height6.3ptdepth-.5ptwidth.8pt\kern-.8pt}
   {\vrule height6.3ptdepth-.5ptwidth.8pt\kern-.8pt}
   {\vrule height4.1ptdepth-.35ptwidth.6pt\kern-.6pt}
   {\vrule height3.1ptdepth-.25ptwidth.5pt\kern-.5pt}}
\def\f@dge{\mathchoice{}{}{\mkern.5mu}{\mkern.8mu}}
\def\b@c#1#2{{\rm \mkern#2mu\vbar\mkern-#2mu#1}}
\def\b@b#1{{\rm I\mkern-3.5mu #1}}
\def\b@a#1#2{{\rm #1\mkern-#2mu\f@dge #1}}
\def\bb#1{{\count4=`#1 \advance\count4by-64 \ifcase\count4\or\b@a A{11.5}\or
   \b@b B\or\b@c C{5}\or\b@b D\or\b@b E\or\b@b F \or\b@c G{5}\or\b@b H\or
   \b@b I\or\b@c J{3}\or\b@b K\or\b@b L \or\b@b M\or\b@b N\or\b@c O{5} \or
   \b@b P\or\b@c Q{5}\or\b@b R\or\b@a S{8}\or\b@a T{10.5}\or\b@c U{5}\or
   \b@a V{12}\or\b@a W{16.5}\or\b@a X{11}\or\b@a Y{11.7}\or\b@a Z{7.5}\fi}}

\catcode`\X=11 \catcode`\@=12




\let\thischap\jobname

\def\partof#1{\csname returnthe#1part\endcsname}
\def\chapof#1{\csname returnthe#1chap\endcsname}

\def\setchapter#1,#2,#3;{%
  \expandafter\def\csname returnthe#1part\endcsname{#2}%
  \expandafter\def\csname returnthe#1chap\endcsname{#3}%
}

\setchapter 300a,A,II.A;
\setchapter 300b,A,II.B;
\setchapter 300c,A,II.C;
\setchapter 300d,A,II.D;
\setchapter 300e,A,II.E;
\setchapter 300f,A,II.F;
\setchapter 300g,A,II.G;
\setchapter  E53,B,N;
\setchapter  88r,B,I;
\setchapter  600,B,III;
\setchapter  705,B,IV;
\setchapter  734,B,V;

\def\cprefix#1{
\edef\theotherpart{\partof{#1}}\edef\theotherchap{\chapof{#1}}%
\ifx\theotherpart\thispart
   \ifx\theotherchap\thischap 
    \else 
     \theotherchap%
    \fi
   \else 
     \theotherchap\fi}

\def\sectioncite[#1]#2{%
     \cprefix{#2}#1}

\edef\thispart{\partof{\thischap}}
\edef\thischap{\chapof{\thischap}}

\def\lastpage of '#1' is #2.{\expandafter\def\csname lastpage#1\endcsname{#2}}


\def\spuriousreset{}


\expandafter\ifx\csname citeadd.tex\endcsname\relax
\expandafter\gdef\csname citeadd.tex\endcsname{}
\else \message{Hey!  Apparently you were trying to
\string\input{citeadd.tex} twice.   This does not make sense.} 
\errmessage{Please edit your file (probably \jobname.tex) and remove
any duplicate ``\string\input'' lines}\endinput\fi

\sectno=-1   
\localtags
\jjtags
\NoBlackBoxes
\define\mr{\medskip\roster}
\define\sn{\smallskip\noindent}
\define\mn{\medskip\noindent}
\define\bn{\bigskip\noindent}
\define\ub{\underbar}
\define\wilog{\text{without loss of generality}}
\define\ermn{\endroster\medskip\noindent}

\define\dbcu{\dsize\bigcup}
\define \nl{\newline}
\magnification=\magstep 1
\documentstyle{amsppt}

{    
\catcode`@11

\ifx\alicetwothousandloaded@\relax
  \endinput\else\global\let\alicetwothousandloaded@\relax\fi

\gdef\subjclass{\let\savedef@\subjclass
 \def\subjclass##1\endsubjclass{\let\subjclass\savedef@
   \toks@{\def\usualspace{{\rm\enspace}}\eightpoint}%
   \toks@@{##1\unskip.}%
   \edef\thesubjclass@{\the\toks@
     \frills@{{\noexpand\rm2000 {\noexpand\it Mathematics Subject
       Classification}.\noexpand\enspace}}%
     \the\toks@@}}%
  \nofrillscheck\subjclass}
} 


\expandafter\ifx\csname alice2jlem.tex\endcsname\relax
  \expandafter\xdef\csname alice2jlem.tex\endcsname{\the\catcode`@}
\else \message{Hey!  Apparently you were trying to
\string\input{alice2jlem.tex}  twice.   This does not make sense.}
\errmessage{Please edit your file (probably \jobname.tex) and remove
any duplicate ``\string\input'' lines}\endinput\fi

\expandafter\ifx\csname bib4plain.tex\endcsname\relax
  \expandafter\gdef\csname bib4plain.tex\endcsname{}
\else \message{Hey!  Apparently you were trying to \string\input
  bib4plain.tex twice.   This does not make sense.}
\errmessage{Please edit your file (probably \jobname.tex) and remove
any duplicate ``\string\input'' lines}\endinput\fi

\def\renewcommand{\newcommand}	       
\edef\cite{\the\catcode`@}%
\catcode`@ = 11
\let\@oldatcatcode = \cite
\chardef\@letter = 11
\chardef\@other = 12
%
%
%
%
\def\@innerdef#1#2{\edef#1{\expandafter\noexpand\csname #2\endcsname}}%
%
%
\@innerdef\@innernewcount{newcount}%
\@innerdef\@innernewdimen{newdimen}%
\@innerdef\@innernewif{newif}%
\@innerdef\@innernewwrite{newwrite}%
%
%
%
\def\@gobble#1{}%
%
%
%
\ifx\inputlineno\@undefined
   \let\@linenumber = \empty 
\else
   \def\@linenumber{\the\inputlineno:\space}%
\fi
%
%
%
\def\@futurenonspacelet#1{\def\cs{#1}%
   \afterassignment\@stepone\let\@nexttoken=
}%
\begingroup 
\def\\{\global\let\@stoken= }%
\\ 
\endgroup
\def\@stepone{\expandafter\futurelet\cs\@steptwo}%
\def\@steptwo{\expandafter\ifx\cs\@stoken\let\@@next=\@stepthree
   \else\let\@@next=\@nexttoken\fi \@@next}%
\def\@stepthree{\afterassignment\@stepone\let\@@next= }%
%
%
%
\def\@getoptionalarg#1{%
   \let\@optionaltemp = #1%
   \let\@optionalnext = \relax
   \@futurenonspacelet\@optionalnext\@bracketcheck
}%
%
%
\def\@bracketcheck{%
   \ifx [\@optionalnext
      \expandafter\@@getoptionalarg
   \else
      \let\@optionalarg = \empty
      \expandafter\@optionaltemp
   \fi
}%
\def\@@getoptionalarg[#1]{%
   \def\@optionalarg{#1}%
   \@optionaltemp
}%
%
%
%
\def\@nnil{\@nil}%
\def\@fornoop#1\@@#2#3{}%
\def\@for#1:=#2\do#3{%
   \edef\@fortmp{#2}%
   \ifx\@fortmp\empty \else
      \expandafter\@forloop#2,\@nil,\@nil\@@#1{#3}%
   \fi
}%
\def\@forloop#1,#2,#3\@@#4#5{\def#4{#1}\ifx #4\@nnil \else
       #5\def#4{#2}\ifx #4\@nnil \else#5\@iforloop #3\@@#4{#5}\fi\fi
}%
\def\@iforloop#1,#2\@@#3#4{\def#3{#1}\ifx #3\@nnil
       \let\@nextwhile=\@fornoop \else
      #4\relax\let\@nextwhile=\@iforloop\fi\@nextwhile#2\@@#3{#4}%
}%
%
%
%
\@innernewif\if@fileexists
\def\@testfileexistence{\@getoptionalarg\@finishtestfileexistence}%
\def\@finishtestfileexistence#1{%
   \begingroup
      \def\extension{#1}%
      \immediate\openin0 =
         \ifx\@optionalarg\empty\jobname\else\@optionalarg\fi
         \ifx\extension\empty \else .#1\fi
         \space
      \ifeof 0
         \global\@fileexistsfalse
      \else
         \global\@fileexiststrue
      \fi
      \immediate\closein0
   \endgroup
}%
%
%
%
%
\def\bibliographystyle#1{%
   \@readauxfile
   \@writeaux{\string\bibstyle{#1}}%
}%
\let\bibstyle = \@gobble
%
%
\let\bblfilebasename = \jobname
\def\bibliography#1{%
   \@readauxfile
   \@writeaux{\string\bibdata{#1}}%
   \@testfileexistence[\bblfilebasename]{bbl}%
   \if@fileexists
      \nobreak
      \@readbblfile
   \fi
}%
\let\bibdata = \@gobble
%
%
\def\nocite#1{%
   \@readauxfile
   \@writeaux{\string\citation{#1}}%
}%
\@innernewif\if@notfirstcitation
%
%
\def\cite{\@getoptionalarg\@cite}%
%
%
\def\@cite#1{%
   \let\@citenotetext = \@optionalarg
   \printcitestart
   \nocite{#1}%
   \@notfirstcitationfalse
   \@for \@citation :=#1\do
   {%
      \expandafter\@onecitation\@citation\@@
   }%
   \ifx\empty\@citenotetext\else
      \printcitenote{\@citenotetext}%
   \fi
   \printcitefinish
}%
\newif\ifweareinprivate
\weareinprivatetrue
\ifx\shlhetal\undefinedcontrolseq\weareinprivatefalse\fi
\ifx\shlhetal\relax\weareinprivatefalse\fi
\def\@onecitation#1\@@{%
   \if@notfirstcitation
      \printbetweencitations
   \fi
   \expandafter \ifx \csname\@citelabel{#1}\endcsname \relax
      \if@citewarning
         \message{\@linenumber Undefined citation `#1'.}%
      \fi
     \ifweareinprivate
      \expandafter\gdef\csname\@citelabel{#1}\endcsname{%
\strut 
\vadjust{\vskip-\dp\strutbox
\vbox to 0pt{\vss\parindent0cm \leftskip=\hsize 
\advance\leftskip3mm
\advance\hsize 4cm\strut\openup-4pt 
\rightskip 0cm plus 1cm minus 0.5cm ?  #1 ?\strut}}
         {\tt
            \escapechar = -1
            \nobreak\hskip0pt\pfeilsw
            \expandafter\string\csname#1\endcsname
             \pfeilso
            \nobreak\hskip0pt
         }%
      }%
     \else  
      \expandafter\gdef\csname\@citelabel{#1}\endcsname{%
            {\tt\expandafter\string\csname#1\endcsname}
      }%
     \fi  
   \fi
   \csname\@citelabel{#1}\endcsname
   \@notfirstcitationtrue
}%
%
%
\def\@citelabel#1{b@#1}%
%
%
\def\@citedef#1#2{\expandafter\gdef\csname\@citelabel{#1}\endcsname{#2}}%
%
%
%
\def\@readbblfile{%
   \ifx\@itemnum\@undefined
      \@innernewcount\@itemnum
   \fi
   \begingroup
      \def\begin##1##2{%
         \setbox0 = \hbox{\biblabelcontents{##2}}%
         \biblabelwidth = \wd0
      }%
      \def\end##1{}
      %
      %
      \@itemnum = 0
      \def\bibitem{\@getoptionalarg\@bibitem}%
      \def\@bibitem{%
         \ifx\@optionalarg\empty
            \expandafter\@numberedbibitem
         \else
            \expandafter\@alphabibitem
         \fi
      }%
      \def\@alphabibitem##1{%
         \expandafter \xdef\csname\@citelabel{##1}\endcsname {\@optionalarg}%
         \ifx\biblabelprecontents\@undefined
            \let\biblabelprecontents = \relax
         \fi
         \ifx\biblabelpostcontents\@undefined
            \let\biblabelpostcontents = \hss
         \fi
         \@finishbibitem{##1}%
      }%
      \def\@numberedbibitem##1{%
         \advance\@itemnum by 1
         \expandafter \xdef\csname\@citelabel{##1}\endcsname{\number\@itemnum}%
         \ifx\biblabelprecontents\@undefined
            \let\biblabelprecontents = \hss
         \fi
         \ifx\biblabelpostcontents\@undefined
            \let\biblabelpostcontents = \relax
         \fi
         \@finishbibitem{##1}%
      }%
      \def\@finishbibitem##1{%
         \biblabelprint{\csname\@citelabel{##1}\endcsname}%
         \@writeaux{\string\@citedef{##1}{\csname\@citelabel{##1}\endcsname}}%
         \ignorespaces
      }%
      %
      %
      \let\em = \bblem
      \let\newblock = \bblnewblock
      \let\sc = \bblsc
      \frenchspacing
      \clubpenalty = 4000 \widowpenalty = 4000
      \tolerance = 10000 \hfuzz = .5pt
      \everypar = {\hangindent = \biblabelwidth
                      \advance\hangindent by \biblabelextraspace}%
      \bblrm
      \parskip = 1.5ex plus .5ex minus .5ex
      \biblabelextraspace = .5em
      \bblhook
      \input \bblfilebasename.bbl
   \endgroup
}%
%
%
\@innernewdimen\biblabelwidth
\@innernewdimen\biblabelextraspace
%
%
%
\def\biblabelprint#1{%
   \noindent
   \hbox to \biblabelwidth{%
      \biblabelprecontents
      \biblabelcontents{#1}%
      \biblabelpostcontents
   }%
   \kern\biblabelextraspace
}%
%
%
%
\def\biblabelcontents#1{{\bblrm [#1]}}%
%
%
\def\bblrm{\rm}%
%
%
\def\bblem{\it}%
%
%
\def\bblsc{\ifx\@scfont\@undefined
              \font\@scfont = cmcsc10
           \fi
           \@scfont
}%
%
%
\def\bblnewblock{\hskip .11em plus .33em minus .07em }%
%
%
\let\bblhook = \empty
%
%
%
\def\printcitestart{[}
\def\printcitefinish{]}
\def\printbetweencitations{, }
\def\printcitenote#1{, #1}
%
%
%
\let\citation = \@gobble
%
%
%
\@innernewcount\@numparams
%
%
\def\newcommand#1{%
   \def\@commandname{#1}%
   \@getoptionalarg\@continuenewcommand
}%
%
%
\def\@continuenewcommand{%
   \@numparams = \ifx\@optionalarg\empty 0\else\@optionalarg \fi \relax
   \@newcommand
}%
%
%
\def\@newcommand#1{%
   \def\@startdef{\expandafter\edef\@commandname}%
   \ifnum\@numparams=0
      \let\@paramdef = \empty
   \else
      \ifnum\@numparams>9
         \errmessage{\the\@numparams\space is too many parameters}%
      \else
         \ifnum\@numparams<0
            \errmessage{\the\@numparams\space is too few parameters}%
         \else
            \edef\@paramdef{%
               \ifcase\@numparams
                  \empty  No arguments.
               \or ####1%
               \or ####1####2%
               \or ####1####2####3%
               \or ####1####2####3####4%
               \or ####1####2####3####4####5%
               \or ####1####2####3####4####5####6%
               \or ####1####2####3####4####5####6####7%
               \or ####1####2####3####4####5####6####7####8%
               \or ####1####2####3####4####5####6####7####8####9%
               \fi
            }%
         \fi
      \fi
   \fi
   \expandafter\@startdef\@paramdef{#1}%
}%
%
%
%
%
\def\@readauxfile{%
   \if@auxfiledone \else 
      \global\@auxfiledonetrue
      \@testfileexistence{aux}%
      \if@fileexists
         \begingroup
            \endlinechar = -1
            \catcode`@ = 11
            \input \jobname.aux
         \endgroup
      \else
         \message{\@undefinedmessage}%
         \global\@citewarningfalse
      \fi
      \immediate\openout\@auxfile = \jobname.aux
   \fi
}%
%
%
\newif\if@auxfiledone
\ifx\noauxfile\@undefined \else \@auxfiledonetrue\fi
%
%
%
%
\@innernewwrite\@auxfile
\def\@writeaux#1{\ifx\noauxfile\@undefined \write\@auxfile{#1}\fi}%
%
%
%
\ifx\@undefinedmessage\@undefined
   \def\@undefinedmessage{No .aux file; I won't give you warnings about
                          undefined citations.}%
\fi
%
%
\@innernewif\if@citewarning
\ifx\noauxfile\@undefined \@citewarningtrue\fi
%
%
%
\catcode`@ = \@oldatcatcode

\def\pfeilso{\leavevmode
            \vrule width 1pt height9pt depth 0pt\relax
           \vrule width 1pt height8.7pt depth 0pt\relax
           \vrule width 1pt height8.3pt depth 0pt\relax
           \vrule width 1pt height8.0pt depth 0pt\relax
           \vrule width 1pt height7.7pt depth 0pt\relax
            \vrule width 1pt height7.3pt depth 0pt\relax
            \vrule width 1pt height7.0pt depth 0pt\relax
            \vrule width 1pt height6.7pt depth 0pt\relax
            \vrule width 1pt height6.3pt depth 0pt\relax
            \vrule width 1pt height6.0pt depth 0pt\relax
            \vrule width 1pt height5.7pt depth 0pt\relax
            \vrule width 1pt height5.3pt depth 0pt\relax
            \vrule width 1pt height5.0pt depth 0pt\relax
            \vrule width 1pt height4.7pt depth 0pt\relax
            \vrule width 1pt height4.3pt depth 0pt\relax
            \vrule width 1pt height4.0pt depth 0pt\relax
            \vrule width 1pt height3.7pt depth 0pt\relax
            \vrule width 1pt height3.3pt depth 0pt\relax
            \vrule width 1pt height3.0pt depth 0pt\relax
            \vrule width 1pt height2.7pt depth 0pt\relax
            \vrule width 1pt height2.3pt depth 0pt\relax
            \vrule width 1pt height2.0pt depth 0pt\relax
            \vrule width 1pt height1.7pt depth 0pt\relax
            \vrule width 1pt height1.3pt depth 0pt\relax
            \vrule width 1pt height1.0pt depth 0pt\relax
            \vrule width 1pt height0.7pt depth 0pt\relax
            \vrule width 1pt height0.3pt depth 0pt\relax}

\def\pfeilsw{ \leavevmode 
            \vrule width 1pt height0.3pt depth 0pt\relax
            \vrule width 1pt height0.7pt depth 0pt\relax
            \vrule width 1pt height1.0pt depth 0pt\relax
            \vrule width 1pt height1.3pt depth 0pt\relax
            \vrule width 1pt height1.7pt depth 0pt\relax
            \vrule width 1pt height2.0pt depth 0pt\relax
            \vrule width 1pt height2.3pt depth 0pt\relax
            \vrule width 1pt height2.7pt depth 0pt\relax
            \vrule width 1pt height3.0pt depth 0pt\relax
            \vrule width 1pt height3.3pt depth 0pt\relax
            \vrule width 1pt height3.7pt depth 0pt\relax
            \vrule width 1pt height4.0pt depth 0pt\relax
            \vrule width 1pt height4.3pt depth 0pt\relax
            \vrule width 1pt height4.7pt depth 0pt\relax
            \vrule width 1pt height5.0pt depth 0pt\relax
            \vrule width 1pt height5.3pt depth 0pt\relax
            \vrule width 1pt height5.7pt depth 0pt\relax
            \vrule width 1pt height6.0pt depth 0pt\relax
            \vrule width 1pt height6.3pt depth 0pt\relax
            \vrule width 1pt height6.7pt depth 0pt\relax
            \vrule width 1pt height7.0pt depth 0pt\relax
            \vrule width 1pt height7.3pt depth 0pt\relax
            \vrule width 1pt height7.7pt depth 0pt\relax
            \vrule width 1pt height8.0pt depth 0pt\relax
            \vrule width 1pt height8.3pt depth 0pt\relax
            \vrule width 1pt height8.7pt depth 0pt\relax
            \vrule width 1pt height9pt depth 0pt\relax
      }


\def\widestnumber#1#2{}

\def\citewarning#1{\ifx\shlhetal\relax 
    \else
    \par{#1}\par
    \fi
}

\def\rm{\fam0 \tenrm}

\def\fakesubhead#1\endsubhead{\bigskip\noindent{\bf#1}\par}



%
%
%

%

\font\textrsfs=rsfs10
\font\scriptrsfs=rsfs7
\font\scriptscriptrsfs=rsfs5

\newfam\rsfsfam
\textfont\rsfsfam=\textrsfs
\scriptfont\rsfsfam=\scriptrsfs
\scriptscriptfont\rsfsfam=\scriptscriptrsfs

\edef\oldcatcodeofat{\the\catcode`\@}
\catcode`\@11

\def\Cal@@#1{\noaccents@ \fam \rsfsfam #1}

\catcode`\@\oldcatcodeofat


\expandafter\ifx \csname margininit\endcsname \relax\else\margininit\fi

\long\def\red#1\endred{}
\long\def\green#1\endgreen{}
\long\def\blue#1\endblue{}
\long\def\private#1\endprivate{}

\def\endred{ \unmatched endred! }
\def\endgreen{ \unmatched endgreen! }
\def\endblue{ \unmatched endblue! }
\def\endprivate{ \unmatched endprivate! }

\ifx\latexcolors\undefinedcs\def\latexcolors{}\fi

\def\emptycs{}
\def\evaluatelatexcolors{%
        \ifx\latexcolors\emptycs\else
        \expandafter\xxevaluate\latexcolors\xxfertig\evaluatelatexcolors\fi}
\def\xxevaluate#1,#2\xxfertig{\setupthiscolor{#1}%
        \def\latexcolors{#2}}


\font\smallfont=cmsl7
\def\rutgerscolor{\ifmmode\else\endgraf\fi\smallfont
\advance\leftskip0.5cm\relax}
\def\setupthiscolor#1{\edef\tmptmpcs{\noexpand\bgroup\noexpand\rutgerscolor
\noexpand\def\noexpand\currentcolor{#1}%
\noexpand}%
\expandafter\let\csname#1\endcsname\tmptmpcs
\def\tmptmpcs{\checkColorUnmatched{#1}\popthecolor}
\expandafter\let\csname end#1\endcsname\tmptmpcs}

\def\checkColorUnmatched#1{\def\expectcolor{#1}%
    \ifx\expectcolor\currentcolor   
    \else \edef\failhere{\noexpand\tryingToClose '\currentcolor' with end\expectcolor}\failhere\fi}

\def\currentcolor{???}

\def\popthecolor{\ifmmode\else\endgraf\fi\egroup}

\expandafter\def\csname#1\endcsname{}

\evaluatelatexcolors

 \let\outerhead\head
 \def\head{\innerhead}
 \let\innerhead\outerhead

 \let\outersubhead\subhead
 \def\subhead{\innersubhead}
 \let\innersubhead\outersubhead

 \let\outersubsubhead\subsubhead
 \def\subsubhead{\innersubsubhead}
 \let\innersubsubhead\outersubsubhead

 \let\outerproclaim\proclaim
 \def\proclaim{\innerproclaim}
 \let\innerproclaim\outerproclaim

 %
 %
 %
 %

\def\demo#1{\medskip\noindent{\it #1.\/}}
\def\enddemo{\smallskip}

\def\remark#1{\medskip\noindent{\it #1.\/}}
\def\endremark{\smallskip}

\pageheight{8.5truein}
\topmatter 
\title{Reflection implies the SCH} \endtitle

\author {Saharon Shelah \thanks {
   The author would like to thank the Israel Science Foundation for
   partial support of this research (Grant no. 242/03). Publication 794. 
\null\newline I would like to thank 
Alice Leonhardt for the beautiful typing.} \endthanks} \endauthor

\affil{Institute of Mathematics\\
 The Hebrew University\\
 Jerusalem, Israel
 \medskip
 Rutgers University\\
 Mathematics Department\\
 New Brunswick, NJ  USA} \endaffil

\keywords  reflection, stationary sets, Singular Cardinal Hypotheses,
 pcf, set theory \endkeywords

\subjclass  03E04, 03E05  \endsubjclass

\abstract  We prove that, e.g., if $\mu > \text{ cf}(\mu) = \aleph_0$
and $\mu > 2^{\aleph_0}$ and every stationary family of countable
subsets of $\mu^+$ reflect in some subset of $\mu^+$ of cardinality
$\aleph_1$ \ub{then} the SCH for $\mu^+$ holds
(moreover, for $\mu^+$, any scale for
$\mu^+$ has a bad stationary set of cofinality $\aleph_1$).
This answers a question of Foreman and Todor\v cevi\'c who 
gets such conclusion from the simultaneous
reflection of four stationary sets. 
\endabstract
\endtopmatter
\document

\newpage

\head {\S0 Introduction} \endhead  \resetall \sectno=0
 \spuriousreset
\bigskip

In \S1 we prove that the strong hypothesis (pp$(\mu) = \mu^+$ for every
singular $\mu$) hence the SCH (singular cardinal hypothesis, that is
$\lambda^\kappa \le \lambda^+ + 2^\kappa$) holds when:
for every $\lambda \ge \aleph_1$ every stationary ${\Cal S} \subseteq
[\lambda]^{\aleph_0}$ reflect in some $A \in [\lambda]^{\aleph_1}$.

This answers a question of Foreman and Todor\v cevi\'c 
\cite{FoTo} where they
proved that the SCH holds for every $\lambda \ge \aleph_1$ when: every
four stationary ${\Cal S}_\ell \subseteq [\lambda]^{\aleph_0},\ell =
1,2,3,4$ reflect simultaneously in some $A \in
[\lambda]^{\aleph_1}$.  They were probably motivated by Velickovi\'c 
\cite{Ve92a} which used another reflection principle: for every stationary
${\Cal A} \subseteq [\lambda]^{\aleph_0}$ there is $A \in
[\lambda]^{\aleph_1}$ such that ${\Cal A} \cap [A]^{\aleph_0}$ contains
a closed unbounded subset, rather than just a stationary set.

The proof here is self-contained modulo two basic
quotations from \cite{Sh:g}, \cite{Sh:f}; we continued \cite{Sh:e},
\cite{Sh:755} in some respects.  
We prove more in \S1.  In particular if 
$\mu > \text{ cf}(\mu) = \aleph_0$ and pp$(\mu) >
\mu^+$ then some ${\Cal A} \subseteq [\mu^+]^{\aleph_0}$ reflect in no
uncountable $A \in [\mu^+]^{\le \mu}$ and see more in the end.

We thank the referee and Shimoni Garti for not few helpful comments.
\mn
For the reader's convenience let us recall some basic definitions.
\definition{\stag{0.1} Definition}  Assume $\theta$ is regular
uncountable (if $\theta = \sigma^+,[B]^{< \sigma} = [B]^{\le \theta}$,
we can use $[B]^{\le \theta}$, the main case is $B = \lambda$)
\mr
\item "{$(a)$}"  ${\Cal A} \subseteq [B]^{< \theta}$ is closed in
$[B]^{< \theta}$, if for every $\{x_\beta:\beta < \alpha\} \subseteq
{\Cal A}$ where $0 < \alpha < \theta$ and $\beta_1 < \beta_2 < \alpha
\Rightarrow x_{\beta_1} \subseteq x_{\beta_2}$, we have $\dbcu_{\beta
< \alpha} x_\beta \in {\Cal A}$
\sn
\item "{$(b)$}"  ${\Cal A}$ is unbounded in $[B]^{< \theta}$, if for
any $y \in [B]^{< \theta}$ we can find $x \in {\Cal A}$, such that $x
\supseteq y$
\sn
\item "{$(c)$}"  ${\Cal A}$ is a club in $[B]^{< \theta}$, if
${\Cal A} \subseteq [B]^{< \theta}$ and (a)+(b) hold for 
${\Cal A}$
\sn
\item "{$(d)$}"  ${\Cal A}$ is stationary in $[B]^{< \theta}$,
or is a stationary subset of $[B]^{< \theta}$ \ub{when} 
${\Cal A} \subseteq [B]^{<\theta}$ and 
${\Cal A} \cap {\Cal C} \ne \emptyset$ for every club
${\Cal C}$ of $[B]^{< \theta}$
\sn
\item "{$(e)$}"  similarly for $[B]^{\le \theta}$ or
$[B]^\theta$ or consider ${\Cal S} \subseteq [B]^{< \theta}$ as 
a subset of $[B]^{\le \theta}$.
\endroster
\enddefinition
\bigskip

\remark{\stag{0.1A} Remark}  Note: if $B = \theta$ then ${\Cal A}
\subseteq [B]^{< \theta}$ is stationary iff ${\Cal A} \cap \theta$ is
a stationary subset of $\theta$.
\endremark
\bigskip

\definition{\stag{0.2} Definition}  Let ${\Cal A} \subseteq [B_1]^{<
\theta}$ and $B_2 \in [B_1]^\mu$.  We say that ${\Cal A}$ reflects in $B_2$
when ${\Cal A} \cap [B_2]^{< \theta}$ is a stationary subset of
$[B']^{< \theta}$.
\enddefinition
\bigskip

\definition{\stag{0.3} Definition}  Let $\kappa$ be a regular
uncountable cardinal, and assume ${\Cal A}$ is a stationary subset of
$[B]^{< \kappa}$.  We define $\diamondsuit_{\Cal A}$ (i.e., the
diamond principle for ${\Cal A}$) as the following assertion:  
\nl
there exists a
sequence $\langle u_a:a \in {\Cal A}\rangle$, such that $u_a \subseteq a$ for
any $a \in {\Cal A}$, and for every $B' \subseteq B$ the set $\{a \in
{\Cal A}:B' \cap a = u_a\}$ is stationary in $[B']^{< \kappa}$.
\enddefinition
\bigskip

\demo{\stag{0.6} Notation}  1) For regular $\lambda > \kappa$ let
$S^\lambda_\kappa = \{\delta < \lambda$: cf$(\delta) = \kappa\}$.
\nl
2) ${\Cal H}(\lambda)$ is the set of $x$ with transitive closure of
   cardinality $< \lambda$.
\nl
3) $<^*_\lambda$ denotes any well ordering of ${\Cal H}(\lambda)$.
\enddemo
\bn
Let us repeat the definition of the next ideal: (see
\cite{Sh:E12}).
\definition{\stag{0.4} Definition}  For $S \subseteq \lambda$ we say
that $S \in \check I[\lambda]$ iff: there is a club $E$ in $\lambda$,
and a sequence $\langle C_\alpha:\alpha < \lambda\rangle$ such that:
\mr
\widestnumber\item{$(iii)$}
\item "{$(i)$}"  $C_\alpha \subseteq \alpha$ for every $\alpha < \lambda$
\sn
\item "{$(ii)$}"  otp$(C_\alpha) < \alpha$
\sn
\item "{$(iii)$}"  $\beta \in C_\alpha \Rightarrow C_\beta = \beta 
\cap C_\alpha$
\sn
\item "{$(iv)$}"  $\alpha \in E \cap S \Rightarrow \alpha =
\sup(C_\alpha)$.
\endroster
\enddefinition
\bigskip

\proclaim{\stag{0.4.3} Claim}  (By \cite{Sh:420} or see \cite{Sh:E12}).
\nl
1) If $\kappa,\lambda$ are regular and $\lambda > \kappa^+$ \ub{then} there
is a stationary $S \subseteq S^\lambda_\kappa$ 
such that $S \in \check I[\lambda]$.
\nl
2) In \scite{0.4} we can add $\alpha \in E \cap S 
\Rightarrow \text{\rm otp}(C_\alpha) = \,\text{\rm cf}(\alpha)$.
\endproclaim
\bigskip

\definition{\stag{0.5} Definition/Observation}  Let ${\Cal A} \subseteq
[\lambda]^\theta$ be stationary and $\lambda \ge \sigma > \theta$ and
$\sigma$ has uncountable cofinality \ub{then} prj$_\sigma({\Cal A}) 
:= \{\sup(a):a \in {\Cal A}\}$, 
and it is a stationary subset of $\sigma$; if $\sigma = \lambda$ we
may omit it.
\enddefinition
\bigskip

\definition{\stag{0.9} Definition}  Let $f_i$ be a function with domain
$\aleph_0$ to the ordinals, for every 
$i \in I$ where $I$ is a set of ordinals. 
We say that the sequence $\bar f = \langle f_i:i \in
I\rangle$ is free, if we can find a sequence $\bar n = \langle n_i:i
\in I\rangle$ of natural numbers such that:  $(i,j \in I) \wedge (i
<j) \wedge (n_i,n_j \le n < \omega) \Rightarrow f_i(n) < f_j(n)$.  We
say that $\bar f$ is $\mu$-free when for every $J \in [I]^{< \mu}$ the
sequence $\bar f \restriction J$ is free.
\enddefinition
\bigskip

\remark{\stag{0.10} Remark}  If we consider ``$\langle f_\alpha:\alpha
\in S\rangle$ for some stationary $S \subseteq \theta$" when $\theta =
\text{ cf}(\theta) > \aleph_0$, then we can assume (without loss of
generality) that $n_i = n(*)$ for every $i \in S$, as we can decrease $S$.
\endremark
\newpage

\head {\S1 Reflection in $[\mu^+]^{\aleph_0}$ and the strong
hypothesis} 
\endhead  \resetall 
 \spuriousreset
\bigskip

\proclaim{\stag{nr.1} The Main Claim}  Assume
\mr
\item "{$(A)$}"   $\lambda = \mu^+$ and $\mu > {\text {\rm cf\/}}(\mu) =
\aleph_0$ and $\aleph_2 \le \mu_* \le \lambda$ (e.g., $\mu_*
= \aleph_2$ which implies that below always $\theta = \aleph_1$) 
\sn
\item "{$(B)$}"  $\bar \lambda = \langle \lambda_n:n < \omega \rangle$ is an
increasing  sequence of regular cardinals $> \aleph_1$ with limit
$\mu$ and $\lambda = { \text{\rm tcf\/}}(\Pi \lambda_n,
<_{J^{\text{bd}}_\omega})$
\sn
\item "{$(C)$}"  $\bar f = \langle f_\alpha:\alpha < \lambda
\rangle$ is $<_{J^{\text{bd}}_\omega}$-increasing cofinal in 
$(\dsize \prod_{n <  \omega} \lambda_n,<_{J^{\text{bd}}_\omega})$
\sn
\item "{$(D)$}"  the sequence $\bar f$ is $\mu_*$-free or at least for 
every cardinal $\theta$ for which 
$\aleph_1 \le \theta = { \text{\rm cf\/}}(\theta) < \mu_*$ the
following is satisfied: if $\theta \le \sigma < \mu_*,{\Cal A} 
\subseteq [\sigma]^{\aleph_0}$ is stationary (recall \scite{0.1A})
and $\langle
\delta_i:i < \theta \rangle$ is an increasing continuous sequence of
ordinals $< \lambda$ \ub{then} for some stationary subfamily ${\Cal
A}_1$ of ${\Cal A}$ (${\Cal A}_1$ is stationary in 
$[\sigma]^{\aleph_0}$ of course)
letting $R_1 = \text{\rm prj}_\theta({\Cal A}_1)$,
see \scite{0.5} we have $\langle f_{\delta_i}:i \in R_1 \rangle$ is
free.  See \scite{0.9} and by \scite{0.10} we can assume that $i \in R_1
\Rightarrow n_i = n(*)$ so $\langle f_{\delta_i}(n):
i \in R_1 \rangle$ is strictly increasing
for every $n \in [n(*),\omega)$.
\ermn
\ub{Then} some stationary ${\Cal A} \subseteq [\lambda]^{\aleph_0}$
does not reflect in any $A \in [\lambda]^{\aleph_1}$ or even in any
uncountable $A \in [\lambda]^{< \mu_*}$, (see Definition \scite{0.2}). 
\endproclaim  
\bigskip

\remark{\stag{nr.1.3} Remark}  0) From the main claim the result on SCH
should be clear from pcf theory (by translating between the pp,cov
and cardinal arithmetic) but we shall give details (i.e. quotes).
\nl
1) Clause (D) from claim \scite{nr.1} is
related to ``the good set of $\bar f$, gd$(\bar f)$ 
contains $S^\lambda_\theta$ modulo the club filter".  
But the clause (D) is stronger.

Note that the good set gd$(\bar f)$ of $\bar f$ is
$\{\delta < \lambda:\aleph_0 < \text{ cf}(\delta) < \mu$ and for
some increasing sequence $\langle \alpha_i:i < \text{
cf}(\delta)\rangle$ of ordinals with limit $\delta$ and sequence $\bar
n = \langle n_i:i < \text{ cf}(\delta)\rangle$ of natural numbers
we have $i < j < \text{ cf}(\delta) \wedge n_i \le n < \omega 
\wedge n_j \le n < \omega \Rightarrow f_{\alpha_i}(n) < 
f_{\alpha_j}(n)$ (so $\langle \cup\{f_{\alpha_i}(n):i
< \text{ cf}(\delta)$ and $n \ge n_i\}:n < \omega\rangle$ is a
$<_{J^{\text{bd}}_\omega}$-eub of $\bar f \restriction \delta)$.

If we use another ideal $J$ say on $\theta < \mu$, the $n_i$ is
replaced by $s_i \in J$.
\nl
2) Recall that by using the silly square (\cite[II,1.5A,pg.51]{Sh:g}), 
if cf$(\mu) \le \theta < \mu,J$ an 
ideal on $\theta$ (e.g. $\theta = \aleph_0,J =   
J^{\text{bd}}_\omega$) and pp$_J(\mu) > \lambda = 
\text{ cf}(\lambda) > \mu$ \ub{then} we can find a sequence $\langle
\lambda_i:i < \theta \rangle$ of regulars $< \mu$ such that $\mu =
\lim_J\langle \lambda_i:i < \theta\rangle$ and tcf$(\dsize \prod_{i < \theta}
\lambda_i,<_J)=\lambda$ and some $\bar f = \langle f_\alpha:\alpha <
\lambda\rangle$ exemplifying it satisfies gd$(\bar f) = \{\delta <
\lambda:\theta < \text{ cf}(\delta) < \mu\}$ and moreover $\bar f$ is
$\mu^+$-free which here means that for every $u \subseteq \lambda$ of
cardinality $\le \mu$ we can find $\langle s_\alpha:\alpha \in
u\rangle$ such that $s_\alpha \in J$, and for $\alpha < \beta$ from
$u$ we have $\varepsilon \in \theta \backslash (s_\alpha \cup s_\beta)
\Rightarrow f_\alpha(\varepsilon) < f_\beta(\varepsilon)$.  This is
stronger than the demand in clause (D).
\nl
3) Also recall that if $\kappa$ is supercompact, $\mu > \kappa >
\theta = \text{ cf}(\mu)$ and $\langle \lambda_i:i < \theta\rangle$ is
an increasing sequence of regulars with limit $\mu,\langle
f_\alpha:\alpha < \lambda \rangle$ exemplifies $\lambda = \mu^+ 
= \text{ tcf}(\dsize \prod_{i < \theta}
\lambda_i,<_{J^{\text{bd}}_\theta})$ then for unboundedly many
$\kappa' \in \kappa \cap \text{ Reg} \backslash
\theta^+$ the set $S^\lambda_{\kappa'} \backslash \text{gd}(\bar f)$ is
stationary.  This is preserved by e.g. Levy$(\aleph_1,< \kappa)$.
\nl
4) For part of the proof (mainly subclaim \scite{nr.4}) we can weaken
clause (D) of the assumption, e.g. in the end demand ``$\Rightarrow
f_{\delta_i}(n) \ne f_{\delta_j}(n)"$ only.  The weakest version of
clause (D) which suffices there is: for any club $C$ of $\theta$ the
set $\cup\{\text{Rang}(f_\alpha):\alpha \in C\}$ has cardinality
$\theta$. 
\endremark
\bn
Before proving \scite{nr.1} we draw some conclusions.
\demo{\stag{nr.2} Conclusion}  1) Assume $\mu > 2^{\aleph_0}$ then
$\mu^{\aleph_0} = \mu^+$ provided that
\mr
\item "{$(A)_\mu$}"  $\mu > \text{ cf}(\mu) = \aleph_0$
\sn
\item "{$(B)_\mu$}"  every stationary ${\Cal A} \subseteq
[\mu^+]^{\aleph_0}$ reflects in some $A \in [\mu^+]^{\aleph_1}$.
\ermn
2) Assume $\lambda \ge \mu_* \ge \aleph_2$.  We can replace $(B)_\mu$ by
\mr
\item "{$(B)_{\mu,\mu_*}$}"   every stationary ${\Cal A} \subseteq
[\mu^+]^{\aleph_0}$ reflects in some uncountable $A \in [\mu^+]^{< \mu_*}$.
\endroster
\enddemo
\bigskip

\demo{Proof}  1) Easily if $\aleph_1 \le \mu' \le \mu$ then
$(B)_{\mu'}$ holds.  Now if $\mu$ is a counterexample, 
\wilog \, $\mu$ is a minimal
counterexample and then by \cite[IX,\S1]{Sh:g} we have 
pp$(\mu) > \mu^+$, hence there is a sequence $\langle
\lambda^0_n:n < \omega \rangle$ of regular cardinals with limit $\mu$ such that
$\mu^{++} = \text{ tcf}(\dsize \prod_{n < \omega}
\lambda^0_n/J^{\text{bd}}_\omega)$, (see \cite{Sh:g}; more
\cite{Sh:E12} or \cite[6.5]{Sh:430}; e.g. using ``no hole for pp" and
the pcf theorem).
Let $\bar f^0 =
\langle f^0_\alpha:\alpha < \mu^{++} \rangle$ witness this.
Hence by \cite[II,1.5A,p.51]{Sh:g} there is $\bar f$ as required in
\scite{nr.1} even a $\mu^+$-free one and 
also the other assumptions there hold so we can
conclude that there exists ${\Cal A} \subseteq [\mu^+]^{\aleph_0}$
which does not reflect in any $A \in [\mu^+]^{\aleph_1}$,
so we get a contradiction to $(B)_\mu$. 
\nl
2) The same proof. \hfill$\square_{\scite{nr.2}}$
\enddemo
\bigskip

\demo{\stag{nr.3} Conclusion}  1) If for every $\lambda > \aleph_1$, every
stationary ${\Cal A} \subseteq [\lambda]^{\aleph_0}$ reflects in some $A \in
[\lambda]^{\aleph_1}$, \ub{then}
\mr
\item "{$(a)$}"  the strong hypothesis (see \cite{Sh:410}, \cite{Sh:420},
\cite{Sh:E12}) holds, i.e. for every singular $\mu$, pp$(\mu) = \mu^+$
and moreover cf$([\mu]^{\text{cf}(\mu)},\subseteq) = \mu^+$ which follows
\sn
\item "{$(b)$}"  the SCH holds.
\ermn
2) Let $\theta \ge \aleph_0$.  We can restrict ourselves to $\lambda >
\theta^+,A \in [\lambda]^{\theta^+}$ (getting the strong
hypothesis and SCH above $\theta$).
\enddemo
\bigskip

\demo{Proof}  1) As in \scite{nr.2}, by \scite{nr.1} we have
$\mu > \text{cf}(\mu) = \aleph_0
\Rightarrow$ pp$(\mu) = \mu^+$, 
this implies clause (a) (i.e. by
\cite[VIII,\S1]{Sh:g}, $\mu > \text{ cf}(\mu) \Rightarrow$ pp$(\mu) =
\mu^+$).  Hence inductively by 
\cite[IX,1.8,pg.369]{Sh:g}, \cite[1.1]{Sh:430} we have 
$\kappa < \mu \Rightarrow \text{ cf}([\mu]^\kappa,\subseteq) 
= \mu$ if cf$(\mu) > \kappa$ and is
$\mu^+$ if $\mu > \kappa \ge \text{ cf}(\mu)$.  This is a consequence
of the strong hypothesis.)  The SCH follows. 
\nl
2) The same proof.    \hfill$\square_{\scite{nr.3}}$
\enddemo
\bigskip

\demo{Proof of \scite{nr.1}}  Let $M^*$ be 
an algebra with universe $\lambda$ and
countably many functions, e.g. all those definable in $({\Cal
H}(\lambda^+),\in,<^*_{\lambda^+},\bar f)$ and are functions
from $\lambda$ to $\lambda$ or just the functions $\alpha \mapsto
f_\alpha(n),\alpha \mapsto \alpha +1$.
\enddemo
\bigskip

\proclaim{\stag{nr.4} Subclaim}  There are $\bar S,S^*,\bar D$ such that:
\mr
\widestnumber\item{$(*)_3(iii)$}
\item "{$(*)_1\,\,\,\,\,\,\,$}"  $\bar S = \langle S_\varepsilon:
\varepsilon < \omega_1 \rangle$ is 
a sequence of pairwise disjoint stationary subsets of $S^\lambda_{\aleph_0}$ 
\sn
\item "{$(*)_2(i)$}"  $S^* \subseteq S^\lambda_{\aleph_1} = \{\delta <
\lambda:{\text{\rm cf\/}}(\delta) = \aleph_1\}$ is stationary and belongs to 
$\check I[\lambda]$
\sn
\item "{$(ii)$}"  if $\delta \in S^*$ \ub{then} there is an increasing
continuous sequence $\langle \alpha_\varepsilon:\varepsilon < \omega_1
\rangle$ of ordinals with limit $\delta$ such that for some sequence $\bar
\zeta = \langle \zeta_\varepsilon:\varepsilon \in R \rangle$ of
ordinals $< \omega_1$, the
set $R \subseteq \omega_1$ is stationary, $\varepsilon \in R
\Rightarrow \alpha_\varepsilon \in S_{\zeta_\varepsilon}$ and $\bar
\zeta$ is with no repetitions
\sn
\item "{$(*)_3(i)$}"  $\bar D = \langle
(D_{1,\varepsilon},D_{2,\varepsilon}):\varepsilon < \omega_1 \rangle$
\sn
\item "{$(ii)$}"  $D_{\ell,\varepsilon}$ is a filter on $\omega$
containing the filter of cobounded subsets of $\omega$
\sn
\item "{$(iii)$}"  if $R_1 \subseteq \omega_1$ is unbounded and $A \in
\cap \{D_{1,\varepsilon}:\varepsilon \in R_1\}$ \ub{then} for some
$\varepsilon \in R_1$ we have $A \ne \emptyset \,\,\text{\rm mod } 
D_{2,\varepsilon}$
\sn
\item "{$(iv)$}"  for each $\varepsilon < \omega_1$ for some $A$ we have $A \in
D_{1,\varepsilon} \and \omega \backslash A \in D_{2,\varepsilon}$.
\endroster
\endproclaim 
\bigskip

\remark{\stag{nr.4A} Remark}  1) For \scite{nr.4} we can assume
(A), (B), (C) of \scite{nr.1} and weaken clause (D): 
because (inside the proof below) necessarily for any stationary
$S^* \subseteq S^\lambda_{\aleph_1}$, which belongs to $\check I[\lambda]$,
we can restrict the demand in (D) of \scite{nr.1} 
for any $\langle \delta_i:i < \theta \rangle$ with
limit in $S^*$.  See more in \cite{Sh:775}. \nl
2) In Subclaim \scite{nr.4} we can demand
$\zeta_\varepsilon = \varepsilon$ in $(*)_2(ii)$.  See the proof.
\nl
3) If we like to demand that each
$D_{\ell,\varepsilon}$ is an ultrafilter (or just have ``$A \in
D_{2,\varepsilon}"$ in the end of $(*)_3(iii)$ of \scite{nr.4}), use
\cite{Sh:E3}.   
\endremark
\bigskip

\demo{Proof of the subclaim \scite{nr.4}}  How do we choose them?

Let $\langle A_i:i \le \omega_1 \rangle$ be a sequence of infinite
pairwise almost disjoint subsets of $\omega$.  Let $D_{1,i} = \{A
\subseteq \omega:A_i \backslash A$ is finite$\},D_{2,i} = \{A
\subseteq \omega:A_j \backslash A$ is finite for all but finitely many
$j < \omega_1\}$, so $D_{2,i}$ does not depend on $i$.  Clearly
$\langle (D_{1,i},D_{2,i}):i < \omega_1\rangle$ satisfies $(*)_3$.

Recall, \scite{0.4.3}, that by \scite{0.4.3} and the fact that
 $\lambda > \aleph_\omega > \aleph_2$, there is a stationary 
$S^* \subseteq S^\lambda_{\aleph_1}$
from $\check I[\lambda]$, and so every stationary $S' \subseteq S^*$ has the
same properties (i.e. is a stationary subset of $\lambda$ which belongs to
$\check I[\lambda]$ and is included in $S^\lambda_{\aleph_1}$).

Let $N \prec ({\Cal H}((2^\lambda)^+),\in,<^*)$ be of cardinality $\mu$
such that $\mu +1 \subseteq N$ and $\{\bar \lambda,\mu,\bar f\}$ belongs
to $N$.  Let $C^* = \cap\{C:C \in N$ is a club of $\lambda\}$, so
clearly $C^*$ is a club of $\lambda$.
For each $h \in {}^\lambda(\omega_1)$ we can try $\bar S^h =
\langle S^h_\gamma:\gamma < \omega_1\rangle$ where $S^h_\gamma = \{\delta <
\lambda:\text{cf}(\delta) = \aleph_0$ and $h(\delta) = \gamma\}$, 
so it is enough to show that
for some $h \in N$, the sequence $\bar S^h$ is as required.  As $\|N\|
< \lambda$, for this it is enough to show that for every 
$\delta \in S^\lambda_{\aleph_1} \cap C^*$ (or just for every $\delta
\in S^* \cap C^*$,
or just for stationarily many $\delta \in S^* \cap C^*$) the demand holds
for $\bar S^h$ for some $h \in ({}^\lambda(\omega_1)) \cap N$.  That
is, $S^{\bar h}$ satisfies $(*)_1$ and $(*)_2(ii)$ of subclaim \scite{nr.4}.
Given any $\delta \in S^\lambda_{\aleph_1} \cap C^*$ let $\langle
\alpha_\varepsilon:\varepsilon < \omega_1 \rangle$ be an increasing
continuous sequence of ordinals with limit $\delta$, \wilog \,
$\varepsilon < \omega_1 \Rightarrow \text{ cf}(\alpha_\varepsilon) =
\aleph_0$,  and by assumption $(D)$ of
\scite{nr.1} for some \footnote{Note that if we require just that
$\langle f_{\alpha_\varepsilon}(n):\varepsilon \in R\rangle$ is
without repetitions, \ub{then} for some stationary subset $R'$ of $R$ the
sequence $\langle f_{\alpha_\varepsilon}(n):\varepsilon \in R'\rangle$
is increasing.} stationary $R \subseteq \omega_1$ and $n = n(*)
< \omega$, the sequence $\langle
f_{\alpha_\varepsilon}(n):\varepsilon \in R \rangle$ is strictly
increasing, so let its limit be $\beta^*$.  So $\beta^* \le \mu$ and
cf$(\beta^*) = \aleph_1$ but $\mu +1 \subseteq N$ hence $\beta^* \in N$.

Note that
\mr
\item "{$(*)_1$}"  for every $\beta' < \beta^*$ the set 
$\{\alpha \in S^\lambda_{\aleph_0}:f_\alpha(n(*)) \in 
[\beta',\beta^*)\}$ is a stationary subset of $\lambda$.
\ermn
[Why?  So assume that $\beta' < \beta^*$ and that the set $S' =
\{\alpha \in S^\lambda_{\aleph_0}:f_\alpha(n(*)) \in 
[\beta',\beta^*)\}$ is not a
stationary subset of $\lambda$.  As $\beta^* +1 \subseteq N$ and $\bar
f \in N$ clearly $S' \in N$ hence there is a club $C'$ of $\lambda$
disjoint to $S'$ which belongs to $N$.  Clearly acc$(C')$ too is a
club of $\lambda$ which belongs to $N$ hence $C^* \subseteq \text{
acc}(C')$ hence $\delta \in \text{ acc}(C')$.  
So $\delta = \sup(C' \cap \delta)$, so $C' \cap
\delta$ is a club of $\delta$.  Recall that $\{\alpha_\varepsilon:\varepsilon
\in R\}$ is a stationary subset of $\delta$ of order type $\aleph_1$.

Now by the choice of $\beta^*$ for some $\varepsilon(*) \in R$ we have
$\beta' \le f_{\alpha_{\varepsilon(*)}}(n(*))$, 
hence $\varepsilon \in R \backslash
\varepsilon(*) \Rightarrow f_{\alpha_\varepsilon}(n(*)) \in
[\beta',\beta^*)$, so $\delta$ has a stationary subset included in
$S'$ hence disjoint to $C'$, contradiction.]
\mr
\item "{$(*)_2$}"  for every $\beta' < \beta^*$ there is $\beta'' \in
(\beta',\beta^*)$ such that $\{\alpha \in S^\lambda_{\aleph_0}:
f_\alpha(n(*)) \in [\beta',\beta'')\}$ is a stationary subset of $\lambda$.
\ermn
[Why?  Follows from $(*)_1$ as $\aleph_1 < \lambda$.]

As $\beta^* \in N$ we can find an increasing continuous sequence $\langle
\beta_\xi:\xi < \omega_1\rangle \in N$ of ordinals with limit
$\beta^*$.  So by $(*)_2$
\mr
\item "{$(*)_3$}"  for every $\xi_1 < \omega_1$ for some $\xi_2 \in
(\xi_1,\omega_1)$ the set $\{\alpha \in
S^\lambda_{\aleph_0}:f_\alpha(n(*)) \in
[\beta_{\xi_1},\beta_{\xi_2})\}$ is stationary.
\ermn
Hence for some unbounded subset $u$ of $\omega_1$ we have
\mr
\item "{$(*)_4$}"  for every $\xi \in u$ the set
$\{\alpha \in S^\lambda_{\aleph_0}:f_\alpha(n(*)) \in 
[\beta_\xi,\beta_{\xi +1})\}$ is a stationary subset of $\lambda$.
\ermn
If $2^{\aleph_1} \le \mu$ then $u = \omega_1$
recalling we demand $\langle \beta_\xi:\xi < \omega_1\rangle \in N$.

We define $h:\lambda \rightarrow \omega_1$ by $h(\alpha) = \zeta$
\ub{iff} for some $\xi \in u$ we have $\zeta = \text{ otp}
(u \cap f_\alpha(n(*)))$ and/or
$\xi = 0 \and f_\alpha(n(*)) \ge \text{ sup}(u)$.
\nl
Clearly $h \in N$ is as required.  So $\bar S = \bar S^h$ as required 
exists.  But maybe $2^{\aleph_1} > \mu$, then after $(*)_3$ we
continue as follows.  Let $\bar C = \langle C_\delta:\delta \in
S^{\aleph_3}_{\aleph_1}\rangle$ be such that $C_\delta$ is a club of
$\delta$ of order type $\omega_1$ which guess clubs, i.e. for every
club $C$ of $\aleph_3$ for stationarily many $\delta \in
S^{\aleph_3}_{\aleph_1}$ we have $C_\delta \subseteq C$, exists by 
\cite[III]{Sh:g}.  Without loss of generality $\bar C \in N$.

Now let $\delta_* \in \text{ acc}(C^*)$ has cofinality $\aleph_3$.
Again $\delta_* \le \mu$ belongs to $N$ hence some increasing
continuous sequence $\langle \alpha_\varepsilon:\varepsilon <
\aleph_3\rangle \in N$ has limit $\delta_*$.  Now for each
$\varepsilon \in S^{\aleph_3}_{\aleph_1}$ we could choose above
$\delta = \alpha_\varepsilon$ hence for some $n_\varepsilon < \omega$
we have $(\forall \beta' < \alpha_\varepsilon)(\exists \beta'' <
\alpha_\varepsilon)[\beta' < \beta'' \wedge
(\exists^{\text{stat}}\gamma \in S^\lambda_{\aleph_0})(\beta' \le
f_\gamma(n_\varepsilon) < \beta'']$.  So for some $n_* < \omega$ the
set $S' := \{\varepsilon < \aleph_3:\text{cf}(\varepsilon) = \aleph_1$
and $n_\varepsilon = n_*\}$ is a stationary subset of
$S^{\aleph_3}_{\aleph_1}$.  It follows that $(\forall \varepsilon <
\aleph_3)(\exists \zeta < \aleph_3)[\varepsilon < \zeta \wedge
(\exists^{\text{stat}} \gamma \in
S^\lambda_{\aleph_0})(\alpha_\varepsilon \le f_\gamma(n_*) <
\alpha_{\varepsilon +1})]$.

Let $\zeta_\varepsilon$ be the minimal $\zeta$ as required above, so
$C = \{\xi < \aleph_3$: if $\varepsilon < \xi$ then $\zeta_\varepsilon
< \xi$ and $\xi$ is a limit ordinal$\}$ is a club of $\aleph_3$.
Hence for some $\varepsilon(*) \in S^{\aleph_3}_{\aleph_1}$ we have
$C_{\varepsilon(*)} \subseteq C$.  Let $u := \{\alpha_\zeta:\zeta \in
C_{\varepsilon(*)}\}$ so clearly $\langle \alpha_\zeta:\zeta \in
u\rangle$ belongs to $N$.  \hfill$\square_{\scite{nr.4}}$
\enddemo
\bigskip

\remark{\stag{nr.4.8} Remark}  Why can't we, in the proof of
\scite{nr.4}, after $(*)_3$, put the instead assuming $2^{\aleph_1}
\le \mu$ use ``as $N \prec ({\Cal H}(2^\lambda)^+),\in,<^*)$ \wilog \, $u
= w_1"$?

The set $u$ chosen above depends on $\delta$, so if $2^{\aleph_1} \le
\mu$ still $u \in N$, but otherwise the ``\wilog \, $u \in N$"
does not seem to be justified.
\endremark
\bigskip

\demo{Continuation of the proof of \scite{nr.1}}
Let $S := \cup\{S_\varepsilon:\varepsilon < \omega_1\}$.  For $\varepsilon <
\omega_1,\delta \in S_\varepsilon$ let 

$$
\align
{\Cal A}^\varepsilon_\delta = \{a:&a \in [\delta]^{\aleph_0} 
\text{ is } M^* \text{-closed, sup}(a) = \delta, \\
  &\text{otp}(a) \le \varepsilon \text{ and}\\
  &(\forall^{D_{1,\varepsilon}} n)(a \cap \lambda_n \subseteq f_\delta(n)) \\ 
  &\text{and } (\forall^{D_{2,\varepsilon}} n)(a \cap \lambda_n \nsubseteq
f_\delta(n))\}
\endalign
$$

\hskip67pt ${\Cal A}^\varepsilon = 
\cup\{{\Cal A}^\varepsilon_\delta:\delta \in S_\varepsilon\}$
\medskip

\hskip67pt ${\Cal A} = \cup\{{\Cal A}^\varepsilon:\varepsilon < \omega_1\}$.
\mn
So

\hskip67pt ${\Cal A} \subseteq [\lambda]^{\aleph_0}$.
\enddemo
\bn
As the case $\mu_* = \aleph_2$ was the original question and its proof is
simpler we first prove it.
\proclaim{\stag{nr.5} Subclaim}  ${\Cal A}$ does not reflect in any $A
\in [\lambda]^{\aleph_1}$.
\endproclaim
\bigskip

\demo{Proof}  So assume $A \in [\lambda]^{\aleph_1}$, let $\langle
a_i:i < \omega_1 \rangle$ be an increasing continuous sequence of
countable subsets of $A$ with union $A$, and let $R = \{i <
\omega_1:a_i \in {\Cal A}\}$, and assume toward contradiction that $R$
is a stationary subset of $\omega_1$.  As every $a \in {\Cal A}$ is
$M^*$-closed, necessarily $A$ is $M^*$-closed
and so \wilog \, each $a_i$ is $M^*$-closed.
\nl
For each $i \in R$ as $a_i \in {\Cal A}$ by the definition of 
${\Cal A}$ we can find $\varepsilon_i
< \omega_1$ and $\delta_i \in S_{\varepsilon_i}$ such that $a_i \in {\Cal
A}^{\varepsilon_i}_{\delta_i}$ hence by the definition of ${\Cal
A}^{\varepsilon_i}_{\delta_i}$ we have otp$(a_i) \le
\varepsilon_i$.  But as $A = \cup\{a_i:i < \omega_1\}$ with $a_i$ countable
increasing with $i$ and  
$|A| = \aleph_1$, clearly for some club $E$ of $\omega_1$
the sequence $\langle \text{otp}(a_i):i \in E \rangle$ is strictly
increasing, hence $i \in E \Rightarrow \text{ otp}(i \cap E) \le
\text{ otp}(a_i)$ so \wilog \, $i \in E \Rightarrow i \le
\text{ otp}(a_i)$ and \wilog \, $i < j \in E \Rightarrow
\varepsilon_i < j \le \text{ otp}(a_j)$.

Now $j \in E \cap R \Rightarrow j \le \text{ otp}(a_j) 
\le \varepsilon_j$ so $\langle \varepsilon_i:i
\in E \cap R \rangle$ is strictly increasing but $\langle
S_\varepsilon:\varepsilon < \omega_1 \rangle$ are pairwise disjoint and
$\delta_i \in S_{\varepsilon_i}$ so $\langle \delta_i:i \in E \cap R
\rangle$ is without repetitions; but $\delta_i = \sup(a_i)$ and for
$i<j$ from $R \cap E$ we have $a_i \subseteq a_j$ which implies that
$\delta_i = \sup(a_i) \le \sup(a_j) = \delta_j$ so
necessarily $\langle \delta_i:i \in R \cap E \rangle$ is strictly
increasing.

As $\sup(a_i) = \delta_i$ for $i \in R \cap E$, clearly $\sup(A) =
\cup\{\delta_i:i \in E \cap R\}$ and let $\beta_i = \text{ Min}(A
\backslash \delta_i)$ for $i < \omega_1$, it is well defined as $\langle
\delta_j:j \in R \cap E\rangle$ is strictly increasing.  Thinning $E$ \wilog 
\mr
\item "{$\circledast_1$}"  $i < j \in E \cap R \Rightarrow \beta_i <
\delta_j \and \beta_i \in a_j$. 
\ermn
Note that, by the choice of $M^*$,
\mr
\item "{$\circledast_2$}"  $i \in E \cap R \wedge
i < j \in E \cap R \Rightarrow \beta_i \in
a_j \Rightarrow \dsize \bigwedge_n (f_{\beta_i}(n) \in a_j) \Rightarrow
\dsize \bigwedge_n (f_{\beta_i}(n) + 1 \in a_j)$.
\ermn
As $\langle \delta_i:i \in E \cap R \rangle$ is (strictly) increasing 
continuous and $R \cap E$ is a stationary subset of $\omega_1$ 
clearly by clause $(D)$ of the assumption of \scite{nr.1}
we can find a stationary $R_1 \subseteq E \cap R$ and $n(*)$ such that $i
\in R_1 \wedge j \in R_1 \wedge i < j \wedge n(*) \le n < \omega
\Rightarrow f_{\delta_i}(n) < f_{\delta_j}(n)$.

Now if $i \in R_1$, let $\bold j(i) =: \text{ Min}(R_1 \backslash (i+1))$, so
$f_{\delta_i} \le_{J^{\text{bd}}_\omega} f_{\beta_i} <_{J^{\text{bd}}_\omega}
f_{\delta_{\bold j(i)}}$ so for some $m_i < \omega$ we have $n \in
[m_i,\omega) \Rightarrow f_{\delta_i}(n) \le f_{\beta_i}(n) <
f_{\delta_{j(i)}}(n)$.  Clearly for some stationary $R_2 \subseteq
R_1$ we have $i,j \in R_2 \Rightarrow m_i = m_j = m(*)$, so possibly
increasing $n(*)$ \wilog \, $n(*) \ge m(*)$; so we have (where Ch$_a
\in \dsize \prod_{n < \omega} \lambda_n$ is defined by Ch$_a(n) = 
\sup(a \cap \lambda_n)$ for any $a \in [\mu]^{< \lambda_0}$):
\mr
\item "{$\circledast_3$}"  for $i<j$ from $R_2$ we have $\bold j(i) \le j$ and
{\roster
\itemitem{ $(\alpha)$ }  $f_{\delta_i} 
\restriction [n(*),\omega) \le f_{\beta_i} \restriction [n(*),\omega)$
\sn
\itemitem{ $(\beta)$ }  $f_{\beta_i} \restriction [n(*),\omega) <
f_{\delta_{\bold j(i)}} \restriction [n(*),\omega) \le
f_{\delta_j} \restriction [n(*),\omega)$
\sn
\itemitem{ $(\gamma)$ }  $f_{\beta_i} \restriction [n(*),\omega) <
\text{ Ch}_{a_j} \restriction [n(*),\omega)$, by $\circledast_2$.
\endroster}
\ermn
Now by the definition of ${\Cal A}^{\varepsilon_i}_{\delta_i}$ as $a_i
\in {\Cal A}^{\varepsilon_i}_{\delta_i} \subseteq 
{\Cal A}^{\varepsilon_i}$ we have
\mr
\item "{$\circledast_4$}"  if $i \in R_2$ then
{\roster
\itemitem{ $(\alpha)$ }  Ch$_{a_i} \le_{D_{1,\varepsilon_i}}
f_{\delta_i}$
\sn
\itemitem{ $(\beta)$ }  $f_{\delta_i} <_{D_{2,\varepsilon_i}} 
\text{ Ch}_{a_i}$.
\endroster}
\ermn
Let $f^* \in \dsize \prod_{n < \omega} \lambda_n$ be $f^*(n) =
\cup\{f_{\beta_i}(n):i \in R_2\}$ if $n \ge n(*)$ and zero otherwise.  
As $f_{\beta_i}(n) \in a_{j(i)}$ for $i \in R_2$ by $\circledast_3(\gamma)$
clearly $n \ge n(*) \Rightarrow f^*(n) \le \sup\{\text{Ch}_{a_i}(n):i \in
R_2\} = \sup(A \cap \lambda_n) = \text{ Ch}_A(n)$ and by
$\circledast_3(\beta)$ we have $n \ge n(*) \Rightarrow 
\text{ cf}(f^*(n)) = \aleph_1$.  Let
$B_1 =: \{n < \omega:n \ge n(*)$ and $f^*(n) = \sup(A \cap
\lambda_n)\}$ and $B_2 =: [n(*),\omega) \backslash B_1$.
As $\alpha \in A \Rightarrow \alpha +1 \in A$ we have $n \in B_1
\Rightarrow A \cap \lambda_n \subseteq f^*(n) = \sup(A \cap
\lambda_n)$.  Also as by the previous sentence $f^* \restriction
[n(*),\omega) \le \text{ Ch}_A \restriction [n(*),\omega)$ clearly
$n \in B_2 \Rightarrow A \cap \lambda_n \nsubseteq f^*(n)$.  
As $\langle a_i:i \in R_2\rangle$ is increasing with union $A$, 
clearly there is $i(*) \in R_2$ such that: 
$n \in B_2 \Rightarrow a_{i(*)} \cap \lambda_n
\nsubseteq f^*(n)$, so as $i \in R_2 \and \alpha \in a_i \Rightarrow
\alpha +1 \in a_i$ we have $i(*) \le i \in R_2 \Rightarrow \text{ Ch}_{a_i}
\restriction B_2 > f_{\delta_i} \restriction B_2$ hence by clause
$\circledast_4(\alpha)$ we have $i \in R_2 \backslash i(*) 
\Rightarrow B_2 = \emptyset$ mod $D_{1,\varepsilon_i} \Rightarrow 
B_1 \in D_{1,\varepsilon_i}$.  
Also by $\circledast_3$ and the choice of $f^*$ and
$B_1$, for each $n \in B_1$ for some club $E_n$ of $\omega_1$ we
have $i \in E_n \cap R_2 \Rightarrow \sup(a_i \cap \lambda_n) =
\sup\{f_{\beta_j}(n):j \in R_2 \cap i\} = \sup\{f_{\delta_j}(n):j \in
R_2 \cap i\} \subseteq f_{\delta_i}(n)$, hence $R_3 = R_2 \cap \bigcap\{E_n
\backslash i(*):n < \omega\}$ is a stationary subset of 
$\omega_1$.  So $n \in
B_1 \and i \in R_3 \Rightarrow a_i \cap \lambda_n \subseteq
f_{\delta_i}(n)$ hence $i \in R_3 \Rightarrow \text{ Ch}_{a_i}
\restriction B_1 \le f_{\delta_i} \restriction B_1$ hence by
$\circledast_4(\beta)$ we have $i \in R_3 \Rightarrow B_1 
= \emptyset$ mod $D_{2,\varepsilon_i}$
hence $i \in R_3 \Rightarrow B_2 \in D_{2,\varepsilon_i}$.

By the choice of $\langle (D_{1,i},D_{2,i}):i < \omega_1 \rangle$ in
\scite{nr.4} as $B_1 \cup B_2$ is a cofinite subset of $\omega,
B_1 \cap B_2 = \emptyset$
(by the choice of $B_1,B_2$, clearly) 
and $R_3 \subseteq \omega_1$ is stationary we get a contradiction, see
$(*)_3(iii)$ of \scite{nr.4}.  \hfill$\square_{\scite{nr.5}}$ 
\enddemo
\bigskip

\proclaim{\stag{nr.6} Subclaim}  ${\Cal A}$ 
is a stationary subset of $[\lambda]^{\aleph_0}$.
\endproclaim
\bigskip

\remark{Remark}  See \cite{RuSh:117}, \cite[XI,3.5,pg.546]{Sh:f},  
\cite[XV,2.6]{Sh:f}.

We give a proof relying only on \cite[XI,3.5,pg.546]{Sh:f}.
In fact, also if we are interested in Ch$_N = \langle \sup(\theta \cap
N):\aleph_0 < \theta \in N \cap \text{ Reg} \rangle,N \prec ({\Cal
H}(\chi),\in)$ we have full control, e.g., if $\bar S = \langle
S_\theta:\aleph_1 \le \theta \in \text{ Reg } \cap \chi \rangle,S_\theta
\subseteq S^\theta_{\aleph_0}$ stationary we can demand
$\aleph_1 \le \theta = \text{ cf}(\theta) \wedge \theta \in N
\Rightarrow \text{ Ch}_N(\theta) \in S_\theta$ and 
control the order of $f^{{\frak a},\lambda}_{\sup(N \cap \lambda)}$ 
and Ch$_N \restriction {\frak a}$.
\endremark
\bigskip

\demo{Proof}  Let $M^{**}$ be an expansion of $M^*$ by countably many
functions; without loss of generality $M^{**}$ has Skolem functions.

Recall that $S^* \subseteq S^\lambda_{\aleph_1}$ is from \scite{nr.4}
so it belongs to $\check I[\lambda]$ and let
$\bar a = \langle a_\alpha:\alpha < \lambda \rangle$ witness it (see
\scite{0.4}, \scite{0.4.3}) so
otp$(a_\alpha) \le \omega_1$ and $\beta \in a_\alpha \Rightarrow a_\beta =
\beta \cap a_\alpha$, and omitting a non-stationary subset of
$S^*$ we have 
$\delta \in S^* \Rightarrow \text{otp}(a_\delta) = \omega_1 \and \delta =
\sup(a_\delta)$.

Let

$$
\align
T^* = \{\eta:&\eta \text{ is a finite sequence of ordinals, } \eta(2n) <
\lambda \text { and} \\
  &\eta(2n+1) < \lambda_n\}.
\endalign
$$
\mn
Let $\lambda_\eta = \lambda$ if $\ell g(\eta)$ is even and
$\lambda_\eta = \lambda_n$ if $\ell g(\eta) = 2n+1$ and let $\bold I_\eta$
be the non-stationary ideal on $\lambda_\eta$ for $\eta \in T^*$, so
$(T^*,\bar{\bold I})$ is well defined where $\bar{\bold I} := \langle \bold
I_\eta:\eta \in T^*\rangle$.

For $\eta \in T^*$, let 
$M_\eta$ be the $M^{**}$-closure of $\{\eta(\ell):\ell < \ell
g(\eta)\}$ so each $M_\eta$ is countable and $\eta \triangleleft \nu
\in T^* \Rightarrow M_\eta \subseteq M_\nu$ and for $\eta \in \text{
lim}(T^*) = \{\eta \in {}^\omega \lambda:\eta \restriction n \in T^*$
for every $n < \omega\}$ let 
$M_\eta = \cup\{M_{\eta \restriction n}:n < \omega\}$,
so it is enough to prove that $M_\eta \in {\Cal A}$ for some $\eta \in
\text{ lim}(T^*)$, more exactly $|M_\eta| \in {\Cal A}$ recall
$M_\eta \subseteq M^{**} \Leftrightarrow M_\eta \prec M^{**}$ as
$M^{**}$ has Skolem functions.  
Let $\bar M = \langle M_\eta:\eta \in T^*\rangle$ and
we can find a subtree $T \subseteq T^*$ such that
\mr
\item "{$\boxtimes$}"   $(T^*,\bold I) \le (T,\bold I)$ and for some
$\varepsilon^* < \omega_1$ we have $\eta \in \text{ lim}(T) \Rightarrow
\text{ otp}(M_\eta) = \varepsilon^*$ \nl
(recalling $(T^*,\bold I) \le (T,\bold I)$ means $T \subseteq
T^*,(\forall \eta \in T^*)(\forall \ell < \ell g(\eta))(\eta
\restriction \ell \in T^*),<> \in T^*$ and $(\forall \eta \in
T)(\{\alpha < \lambda_\eta:\eta \char 94 \langle \alpha\rangle \in
T^*\} \ne \emptyset$ mod $\bold I_\eta$, i.e. is stationary)). 
\ermn
Why?  As lim$(T^*) = \cup\{\bold B_\varepsilon:\varepsilon <
\omega_1\}$, see $\boxtimes_4$ below, and by $\boxtimes_1$ below each $\bold
B_\varepsilon$ is a Borel subset of lim$(T^*)$ and note that
$\boxtimes$ says the $(\exists \varepsilon < \omega_1)(\exists T)[(T^*,\bold I)
\le (T,\bold I) \cap \text{ lim}(T) \subseteq \bold B_\varepsilon)$.
The existence of such $\varepsilon$ is, e.g.,  
\cite[XI;3.5,p.546]{Sh:f}; the reader may ask to
justify the sets being Borel, so let $u_\eta$ be the universe of
$M_\eta$, a countable set of ordinals.

So we use
\mr
\item "{$\boxtimes_1$}"  for any $\varepsilon < \omega_1$ the set
$\bold B_\varepsilon = \{\eta \in \lim(T):\text{otp}(u_\eta) 
= \varepsilon\}$ is a Borel set.
\ermn
[Why?  Without loss of generality $u_\eta \ne \emptyset$ and let
$\langle \alpha_{\eta,n}:n < \omega\rangle$ enumerate the members of
$u_\eta$ and for $n_1,n_2 < \omega$ and $m_1,m_2 <\omega$ let 
$\bold B_{n_1,n_2,m_1,m_2} := \{\eta \in \lim(T^*):\alpha_{\eta
\restriction n_1,m_1} < \alpha_{\eta \restriction n_2,m_2}\}$.

Clearly
\mr
\item "{$\boxtimes_2$}"  $\bold B_{n_1,n_2,m_1,m_2}$ is an open subset of
lim$(T^*)$
\sn
\item "{$\boxtimes_3$}"  there is a $\Bbb L_{\omega_1,\omega}$ sentence
$\psi_\varepsilon$ in the vocabulary consisting of
$\{p_{n_1,n_2,m_1,m_2,\ell}:n_1,n_2,m_1,m_2 < \omega\}$
such that: the $p$'s are propositional variables 
(i.e. $0$-place predicates) and
if $\langle \alpha_{n,m}:n,m < \omega\rangle$ is a sequence
of ordinals and $p_{n_1,n_2,m_1,m_2}$ is assigned the truth value of
$\alpha_{n_1,m_1} < \alpha_{n_2,m_2}$ then $\gamma = \text{
otp}\{\alpha_{n,m}:n,m <\omega\}$ iff $\psi_\varepsilon$ is assigned the truth
value true
\sn
\item "{$\boxtimes_4$}"  lim$(T^*) = \cup\{\bold B_\varepsilon:
\varepsilon < \omega_1\}$.
\ermn
[Why?  As otp$(M_\eta \cap \lambda) < \|M_\eta\|^+ = \aleph_1$.
Together $\boxtimes$ should be clear.]

Note that for every $\eta \in T^*$ of length $2n+2$ we have $\eta
\trianglelefteq \nu \in T^* \Rightarrow \bold I_\nu$ is $\lambda^+_n$-complete.
As we can shrink $T$ further by \cite[XI,3.5,pg.346]{Sh:f} without loss
of generality
\mr
\item "{$\circledast$}"  for every $n < \omega$ and $\eta \in T \cap
{}^{2n+2}\lambda$ for some $\alpha = \alpha_\eta < \lambda_n$ we have:
if $\eta \triangleleft \nu \in \lim(T)$ then $\alpha_\eta =
\sup(\lambda_n \cap M_\nu)$.
\ermn
[Why?  As above applied to each $T' = \{\rho \in 
{}^{\omega >}\lambda:\eta \char 94 \rho \in T\}$.]

Let $\chi = (2^\lambda)^+$ and  
$N^*_\alpha \prec {\frak B} = ({\Cal H}(\chi),\in,<^*_\chi)$ for
$\alpha < \lambda$ be
increasing continuous, $\|N^*_\alpha\| = \mu,\alpha \subseteq
N^*_\alpha,\langle N^*_\beta:\beta \le \alpha \rangle \in N^*_{\alpha
+1}$ and $(T,\bold I,\bar M,\bar a,\bar f,\bar \lambda,\mu) \in
N^*_\alpha$, clearly possible and $E = \{\delta < \lambda:N^*_\delta
\cap \lambda = \delta\}$ is a club of $\lambda$, hence we can find
$\delta(*) \in S^* \cap E$, so $a_{\delta(*)}$ is well defined.  
Let $\bar N^* = \langle N^*_\alpha:\alpha < \lambda\rangle$.
Let
$C_{\delta(*)}$ be the closure of $a_{\delta(*)}$ as a subset of
$\delta(*)$ in the order topology and 
let $\langle \alpha_\varepsilon:\varepsilon < \omega_1
\rangle$ list $C_{\delta(*)}$ in increasing order, so is increasing continuous.

We define $N_\varepsilon$ by induction on $\varepsilon < \omega_1$ by:
\mr
\item "{$(*)_0$}"  $N_\varepsilon$ is the Skolem hull in ${\frak B}$
of
\nl
$\{\alpha_\zeta:\zeta < \varepsilon\} \cup \{\langle N_\xi:\xi < \zeta
\rangle,\bar N^* \restriction \zeta:\zeta < \varepsilon\} \cup 
\{(T,\bold I,\bar M,\bar a,\bar f,\bar \lambda,\mu)\}$.
\ermn
Let
\mr
\item "{$(*)_1$}"  $g_\varepsilon \in \dsize \prod_{n < \omega}
\lambda_n$ be defined by $g_\varepsilon(n) = \text{ sup}(N_\varepsilon
\cap \lambda_n)$.
\ermn
Clearly
\mr
\item "{$(*)_2$}"   $(a) \quad \langle N_\zeta:\zeta \le \varepsilon\rangle \in
N^*_{\delta(*)}$ and even $\in N_\xi$ for every $\xi \in [\varepsilon
+1,\omega_1)$
\sn
\item "{${{}}$}"  $(b) \quad C_{\delta(*)} \cap (\alpha_\varepsilon +1)$ and
$a_{\alpha_\varepsilon}$ belongs to $N_\xi$ for $\xi \in [\varepsilon
+1,\omega_1)$. 
\ermn
[Why?  For clause (a), $\langle N_\zeta:\zeta \le \varepsilon\rangle$
appear in the set whose Skolem Hull is $N_\xi$.  For clause (b)
because $\bar a \in N^*_{\delta(*)}$ and $\alpha \in a_{\delta(*)}
\Rightarrow a_\alpha = a_{\delta(*)} \cap \alpha$ and $C_{\delta(*)}
\cap (\alpha_\varepsilon +1) =$ the closure of $a_{\alpha_{\varepsilon
+1}} \cap (\alpha_\varepsilon +1)$.]
\mn
Let $e = \{\varepsilon < \omega_1:\varepsilon$ is a limit ordinal and
$N_\varepsilon \cap \omega_1 = \varepsilon\}$.  So
\mr
\item "{$(*)_3$}"  $(a) \quad e$ is a club of $\omega_1$,
\sn
\item "{${{}}$}"  $(b) \quad$ if $\varepsilon \in e$ then
sup$(N_\varepsilon \cap \lambda) = \alpha_\varepsilon =
N^*_{\alpha_\varepsilon} \cap \lambda,N_\varepsilon \subseteq
N^*_{\alpha_\varepsilon}$ and $\varepsilon < \zeta < \omega_1 \Rightarrow
N_\varepsilon \in N_\zeta$
\ermn
hence
\mr
\item "{$(*)_4$}"  if $\varepsilon +2 < \zeta \in e$ then
$g_\varepsilon,g_{\varepsilon +1} \in N_{\varepsilon +2} \prec N_\zeta$.
\ermn
Now $\bar f$ is increasing and cofinal in $\dsize \prod_{n < \omega} \lambda_n$
hence
\mr
\item "{$(*)_5$}"   if $\varepsilon < \zeta \in e$ then $g_\varepsilon
<_{J^{\text{bd}}_\omega} f_{\alpha_\zeta}$ and $f_{\alpha_\varepsilon}
<_{J^{\text{bd}}_\omega} g_\zeta$.
\ermn
Also easily
\mr
\item "{$(*)_6$}"  if $\varepsilon < \zeta \in e$ then $g_\varepsilon
< g_\zeta$.
\ermn
For $n < \omega,\varepsilon < \omega_1$ let $N_{\varepsilon,n+1}$ be
the Skolem hull inside ${\frak B}$ of $N_\varepsilon \cup \lambda_n$
and let $N_{\varepsilon,0} = N_\varepsilon$.  Easily
\mr
\item "{$(*)_7$}"  if $n \le m < \omega$ and $\varepsilon < \omega_1$ then
$g_\varepsilon(m) =\sup(N_{\varepsilon,n} \cap \lambda_m)$.
\ermn
Recall that $\varepsilon^*$ is the order type of $M_\eta \cap \lambda$
for every $\eta \in \text{ lim}(T)$.
\nl
Choose $\varepsilon \in \text{ acc}(e)$ such that $\varepsilon >
\varepsilon^*,\alpha_\varepsilon \in S_\zeta$ for some $\zeta \in
[\varepsilon,\omega_1)$ (possible by subclaim \scite{nr.4} particularly
clause $(*)_2(ii)$) and choose $\varepsilon_k \in
e \cap \varepsilon$ for $k < \omega$ such that $\varepsilon_k <
\varepsilon_{k+1} < \varepsilon = \cup\{\varepsilon_\ell:\ell < \omega\}$.
We also choose $n_k$ by induction on $k < \omega$ such that
\mr
\item "{$(*)_8$}"  $(a) \quad n_\ell < n_k < \omega$ for $\ell < k$
\sn
\item "{${{}}$}"   $(b) \quad g_{\varepsilon_{k+1}} 
\restriction [n_k,\omega) < f_{\alpha_\varepsilon} \restriction [n_k,\omega)$.
\ermn
[Why is this choice possible?  By $(*)_5$.]
\mn
Stipulate $n_{-1} = 0$.

Let $B_1 \in D_{1,\zeta}$ be such that $B_2 = \omega \backslash B_1
\in D_{2,\zeta}$, exists by clause $(*)_3(iv)$ of subclaim \scite{nr.4}.

Now we choose $\eta_n$ by induction on $n < \omega$ such that
\mr
\item "{$\boxdot$}"  $(a) \quad \eta_n \in T$ and $\ell g(\eta_n) = n$
\sn
\item "{${{}}$}"  $(b) \quad m < n \Rightarrow \eta_m \triangleleft
\eta_n$
\sn
\item "{${{}}$}"  $(c) \quad$ if $n \in [n_{k-1},n_k)$ then
$\eta_{2n},\eta_{2n+1} \in N_{\varepsilon_k,n}$
\sn
\item "{${{}}$}"  $(d) \quad$ if $n \in [n_{k-1},n_k)$ then
$\eta_{2n+1}(2n) = \text{ Min}\{\alpha < \lambda:\eta_{2n} \char 94 \langle
\alpha \rangle \in T$ and 
\nl

\hskip25pt $\alpha \ge \alpha_{\varepsilon_{k-1}}$ if $k>0\}$
\sn
\item "{${{}}$}"  $(e) \quad$ if $n \in [n_{k-1},n_k)$ and $n \in B_1$
then $\eta_{2n+2}(2n+1) = \text{ Min}\{\alpha < \lambda_n$:
\nl

\hskip25pt $\eta_{2n+1} \char 94 \langle \alpha \rangle \in T\}$
\sn
\item "{${{}}$}"  $(f) \quad$ if $n \in [n_{k-1},n_k)$ and $n \in B_2$
then $\eta_{2n+2}(2n+1) = \text{ Min}\{\alpha < \lambda_n$:
\nl

\hskip25pt $\eta_{2n+1} \char 94 \langle \alpha \rangle \in T$ and $\alpha >
f_{\alpha_\varepsilon}(n)\}$.
\ermn
No problem to carry the induction.
\mn
[Clearly if $\eta_n$ is well defined then $\eta_{n+1}(n)$ is well
defined (by clause (c) or (d) or (e) according to the case; hence
$\eta_{n+1} \in T \cap {}^{n+1}\lambda$ is well defined by why clause
(c) holds, i.e. assume $n \in [n_{k-1},n_k)$, why
$\eta_{2n},\eta_{2n+1} \in N_{\varepsilon_{k,n}}$?
\bn
\ub{Case 1}:  If $n=0$, then $\eta_{2n} = <> \in N_{\varepsilon_k,n}$
trivially.
\bn
\ub{Case 2}:  $\eta_{2n}$ is O.K. hence $\in N_{\varepsilon_{k,n}}$
and show $\eta_{2n+1} \in N_{\varepsilon_{k,n}}$.
\nl
[Why?  Because $N_{\varepsilon_k,n} \prec {\frak B}$, if $k=0$ as
$\eta_{2n+2}(2)$ is defined from $\eta_{2n}$ and $T$ both of which
belongs to $N_{\varepsilon_k,n}$.  If $k>0$ we have to check that also
$\alpha_{\varepsilon_{k-1}} \in N_{\varepsilon_k,n}$ which holds by $(*)_0$.
\bn
\ub{Case 3}:  $\eta_{2n+1}$ is O.K. so $\in N_{\varepsilon_k,n}$
and we have to show $\eta_{2n+2} \in N_{\varepsilon_k,n+1}$.

As $\eta_{2n+2}(n) < \lambda_n \subseteq N_{\varepsilon_k,n+1}$ this
should be clear.]  

Let $\eta = \cup\{\eta_n:n < \omega\}$.  Clearly $\eta \in \text{
lim}(T)$ hence $u =: |M_\eta| \in [\lambda]^{\aleph_0}$ and $M_\eta
\subseteq M^{**}$, hence it is enough to prove that $u \in {\Cal A}$.

Now
\mr
\item "{$\circledast_1$}"  sup$(u) \le \alpha_\varepsilon$
\nl
[Why?  As $\eta_n$ belongs to the Skolem hull of $N_\varepsilon \cup
\mu \subseteq N^*_{\alpha_\varepsilon}$ hence $M_{\eta_n}
\subseteq N_\varepsilon \subseteq N^*_{\alpha_\varepsilon}$ 
and $N^*_{\alpha_\varepsilon} \cap
\lambda = \alpha_\varepsilon$ as $\alpha_\varepsilon \in E$.]
\sn
\item "{$\circledast_2$}"  sup$(u) \ge \alpha_{\varepsilon_n}$, for
every $n < \omega$
\nl
[by clause (d) of $\boxdot$]
\sn
\item "{$\circledast_3$}"  sup$(u) = \alpha_\varepsilon$
\nl
[Why?  By $\circledast_1 + \circledast_2$]
\sn
\item "{$\circledast_4$}"  $\alpha_\varepsilon \in S_\zeta$ and $\zeta
\ge \varepsilon > \varepsilon^* = \text{otp}(u)$
\nl
[Why?  By the choice of $\varepsilon$]
\sn
\item "{$\circledast_5$}"  if $n \ge n_0,n > 0$ and $n \in B_1$ then $u
\cap \lambda_n \subseteq f_{\alpha_\varepsilon}(n)$
\nl
[Why?  By the choice of $\eta_{2n+2}(2n+1)$, i.e., let $k$ be such
that $n \in [n_{k-1},n_k)$, so $\eta_{2n+1} \in N_{\varepsilon_k,n}$
by clause (c) and by clause (e) of $\boxdot$ we have 
$\eta_{2n+2}(2n+1) \in \lambda_n \cap 
N_{\varepsilon_k,n}$ hence by $\otimes$ above, as $\eta \in \text{
lim}(T)$ we have $\alpha_{\eta \restriction (2n+2)} =
\alpha_{\eta_{2n+2}} = \sup(u \cap \lambda_n)$ and as $\bar M \in
N_{\varepsilon_k,n}$ we have $\alpha_{\eta_{2n+2}} \in
N_{\varepsilon_k,n}$ so sup$(u \cap \lambda_n) = \alpha_{\eta_{2n+2}}  
 < \sup(N_{\varepsilon_k,n} \cap \lambda_n)$ but the latter is equal to
sup$(N_{\varepsilon_k} \cap \lambda_n)$ by $(*)_7$ which is equal to 
$g_{\varepsilon_k}(n)$
which is $< f_{\alpha_\varepsilon}(n)$ by $(*)_8$, as required.]
\sn
\item "{$\circledast_6$}"  if $n \ge n_1$ and $n \in B_2$ then $u \cap
\lambda_n \nsubseteq f_{\alpha_\varepsilon}(n)$
\nl
[Why?  By the choice of $\eta_{2n+2}(2n+1)$.]
\ermn
So we are done.  \hfill$\square_{\scite{nr.6}}$
\enddemo
\bn
This (i.e., \scite{nr.5} + \scite{nr.6}) is 
enough for proving \scite{nr.1} in the case $\mu_* = \aleph_2$.  In
general we should replace \scite{nr.5} by the following claim.
\proclaim{\stag{nr.7} Claim}  The family ${\Cal A}$ does not reflect in any
uncountable $A \in [\lambda]^{< \mu_*}$.
\endproclaim
\bigskip

\demo{Proof}  Assume $A$ is a counterexample.

Trivially 
\mr
\item "{$\circledast_0$}"  $A$ is $M^*$-closed.
\ermn
For $a \in {\Cal A}$ let $(\delta(a),\varepsilon(a))$ be such that $a \in
{\Cal A}^{\varepsilon(a)}_{\delta(a)}$ hence $\delta(a) = \sup(a)$, otp$(a)
\le \varepsilon(a)$.  Let ${\Cal A}^- = {\Cal A} \cap
[A]^{\aleph_0}$ and let $\Gamma = \{\delta(a):a \in {\Cal A}^-\}$.
Of course, $\Gamma \ne \emptyset$.  Assume that $\delta_n \in \Gamma$
for $n < \omega$ so let $\delta_n = \delta(a_n)$ where $a_n \in {\Cal A}$ 
so necessarily $\delta_n \in S_{\varepsilon(a_n)}$.  As $A$ is uncountable
we can find a countable $b$ such that $a_n \subseteq b \subseteq A$ and
$\varepsilon(a_n) < \text{ otp}(b)$ for every $n < \omega$
and as ${\Cal A}^- \subseteq [A]^{\aleph_0}$ is stationary 
we can find $c$ such that $b
\subseteq c \in {\Cal A}^-$; so $\varepsilon(c) \ge \text{ otp}(c) \ge \text{
otp}(b) > \varepsilon(a_n) \and \delta_n \in S_{\varepsilon(a_n)} 
\and \delta(a_n) =
\delta_n \le \sup(a_n) \le \sup(c) = \delta(c)$ for each $n < \omega$.  So if
$\delta(a_n) = \delta_n = \delta(c),n < \omega$ necessarily $\varepsilon(a_n) 
= \varepsilon(c)$ contradiction so $\delta_n \ne \delta(c)$; hence
$\delta(c) > \delta(a_n)$ and, of course, $\delta(c) \in \Gamma$ so $n
< \omega \Rightarrow \delta_n < \delta(c) \in 
\Gamma$.  As $\delta_n$ for $n < \omega$ were
any members of $\Gamma$, clearly $\Gamma$ has no
last element, and let $\delta^* = \sup(\Gamma)$.  Similarly cf$(\delta^*)
= \aleph_0$ is impossible, so clearly cf$(\delta^*) > \aleph_0$ and
let $\theta = \text{ cf}(\delta^*)$ so $\theta \le |A| < \mu_*$ and
$\theta$ is a regular uncountable cardinal. \nl
As $a \in {\Cal A}^\varepsilon_\delta \Rightarrow \sup(a) = \delta$ and
${\Cal A}^- \subseteq [A]^{\aleph_0}$ is stationary 
clearly $A \subseteq \delta^* = \sup(A) = \sup(\Gamma)$.  
Let $\langle \delta_i:i < \theta \rangle$ 
be increasing continuous with limit $\delta^*$ and if
$\delta_i \in S_\varepsilon$ then we let $\varepsilon_i = \varepsilon$.

For $i < \theta$ let $\beta_i = \text{ Min}(A \backslash \delta_i)$,
so $\delta_i \le \beta_i < \delta^*,\beta_i \in A$ and $i < j < \theta
\Rightarrow \beta_i \le \beta_j$.
But $i < \theta \Rightarrow \beta_i < \delta^* \Rightarrow (\exists j)(i
< j < \theta \wedge \beta_i < \delta_j)$ so 
for some club $E_0$ of $\theta$ we
have $i < j \in E_0 \Rightarrow \beta_i < \delta_j \le \beta_j$; as we
can replace $\langle \delta_i:i < \theta \rangle$ by $\langle \delta_i:i
\in E_0 \rangle$ \wilog \, $\beta_i < \delta_{i+1}$ hence $\langle
\beta_i:i < \theta \rangle$ is strictly increasing. 
\nl
Let $A^- := \{\beta_i:i < \theta\}$ and let $H:[\theta]^{\aleph_0}
\rightarrow \theta$ be $H(b) = \text{ sup}\{i:\beta_i \in b\}$ and let
$J := \{R \subseteq \theta$: the family $\{b \in {\Cal A}^-:H(b) \in
R\} = \{b \in {\Cal A}^-:\sup(\{i < \theta:\beta_i \in b\}) \in
R\}$ is not a stationary subset of $[A^-]^{\aleph_0}\}$. 
\nl
Clearly
\mr
\item "{$\circledast_1$}"   $J$ is an $\aleph_1$-complete ideal on $\theta$
extending the non-stationary ideal and $\theta \notin J$ by the
definition of the ideal
\sn
\item "{$\circledast_2$}"   if $B \in J^+$ (i.e., $B \in {\Cal
P}(\theta) \backslash J$) \ub{then} $\{a \in {\Cal A}^-:H(a) \in B\}$
is a stationary subset of $[\theta]^{\aleph_0}$.
\ermn
By clause (D) of the assumption of \scite{nr.1}, for some stationary
$R_1 \in J^+$ and $n_i < \omega$ for $i \in R_1$ we have 
\mr
\item "{$\circledast_3$}"   if $i < j$ are from $R_1$
and $n \ge n_i,n_j$ (but $n < \omega$) then $f_{\beta_i}(n) <
f_{\beta_j}(n)$.
\ermn
Recall that 
\mr
\item "{$\circledast_4$}"  $i < j \in R_1 \Rightarrow \beta_i <
\delta_j$. 
\ermn
Now if $i \in R_1$, let $j(i) = \text{ Min}(R_1 \backslash (i+1))$, so
$f_{\delta_i} \le_{J^{\text{bd}}_\omega} f_{\beta_i} <_{J^{\text{bd}}_\omega}
f_{\delta_{\bold j(i)}}$ hence for some $m_i < \omega$ we have $n \in
[m_i,\omega) \Rightarrow f_{\delta_i}(n) \le f_{\beta_i}(n) <
f_{\delta_{\bold j(i)}}(n)$.  Clearly for some $n(*)$ satisfying
$\lambda_{n(*)} > \theta$ and $R_2 \subseteq R_1$ from $J^+$ we have 
$i \in R_2 \Rightarrow n_i,m_i \le n(*)$, so
\mr
\item "{$\circledast_5$}"  for $i<j$ in $R_2$ we have
{\roster
\itemitem{ $(\alpha)$ }  $f_{\delta_i} 
\restriction [n(*),\omega) \le f_{\beta_i}[n(*),\omega)$
\sn
\itemitem{ $(\beta)$ }  $f_{\beta_i} \restriction [n(*),\omega) <
f_{\delta_j} \restriction [n(*),\omega)$.
\endroster}
\ermn
Let $f^* \in \dsize \prod_{n < \omega} \lambda_n$ be defined by
$f^*(n) = \cup \{f_{\delta_i}(n):i \in R_2\}$ if $n \ge n(*)$ and zero
otherwise.  Clearly $f^*(n) \le \sup(A \cap \lambda_n)$ for $n <
\omega$.

Let ${\Cal A}' = \{a \in {\Cal A}^-:(\forall i < \theta)(i \in a
\equiv \beta_i \in a \equiv \delta_i \in \beta_i)
\sup\{i \in R_2:\beta_i \in a\} =
\sup\{i:\beta_i \in a\} = \sup(a \cap \theta) \in R_2$ and 
sup$(A \cap \lambda_n) > f^*(n)
\Rightarrow a \cap \lambda_n \nsubseteq f^*(n)\}$.  As $R_2 \in J^+$ 
clearly ${\Cal A}'$ is a stationary subset of $[A]^{\aleph_0}$. 
\nl
Let $R_3 = \{i \in R_2:i = \sup(i \cap R_2)\}$ so $R_3 \subseteq
R_2,R_2 \backslash R_3$ is a non-stationary subset of $\theta$ (hence
belongs to $J$) and $a \in {\Cal A}' \Rightarrow \sup(a) 
\in \{\delta_i:i \in R_3\}$. \nl
Let

$$
\align
{\Cal A}^* = \biggl\{ a \in [A]^{\aleph_0}:&(a) \quad \beta_{\text{min}(R_2)}
\in a \text{ and } a \text{ is } M^*\text{-closed} \\
  &(b) \quad \text{if } i \in R_2 \and j =
\text{ Min}(R_2 \backslash (i+1)) \text{ \ub{then} } 
[a \nsubseteq \delta_i \Rightarrow a \nsubseteq \delta_j] \text{ and} \\
 &\qquad \,n \in [n(*),\omega) \and a \cap \lambda_n \nsubseteq
f_{\delta_i}(n) \Rightarrow a \cap 
\lambda_n \backslash f_{\delta_j}(n) \ne \emptyset \\
  &(c) \quad \text{if } i < \theta \and n \in 
[n(*),\omega) \text{ then } (\exists \gamma)(\beta_i \le \gamma \in a)
\equiv \\
  &\qquad \, (\exists j)(i < j < \theta
\and \beta_j \in a) \equiv (\exists \gamma)(f_{\beta_i}(n) \le \gamma
\in a \cap f^*(n)) \text{ and} \\
  &(d) \quad \text{if } A \cap \lambda_n \nsubseteq f^*(n) \text{ then
} a \cap \lambda_n \nsubseteq f^*(n) \\
  &\qquad \, \text{ but } (\forall \gamma \in a)(\gamma +1 \in a) \\
  &\qquad \, \text{hence sup } (a \cap \lambda_n) > f^*(n) \biggr\}.
\endalign
$$
\mn
Clearly ${\Cal A}^*$ is a club of 
$[A]^{\aleph_0}$ (recall that $A$ is $M^*$-closed).  But if $a
\in {\Cal A}^* \cap {\Cal A}'$, then for some limit ordinal $i \in R_3
\subseteq  \theta$ we have $a \subseteq \sup(a) = \delta_i$ and 
$n \in [n(*),\omega) \Rightarrow \sup(a \cap f^*(n)) = \sup(a \cap \cup
\{f_{\delta_j}(n):j \in R_2\})$.
\nl
Let 

$$
B_1 = \{n:n(*) \le n < \omega \text{ and } A \cap \lambda_n \subseteq
f^*(n) = \sup(A \cap \lambda_n)\}.
$$

$$
B_2 = \{n:n(*) \le n < \omega \text{ and } f^*(n) < \sup(A \cap
\lambda_n)\}.
$$
\mn
Clearly $B_1,B_2$ are disjoint with union $[n(*),\omega)$ 
recalling $\alpha \in A \Rightarrow \alpha + 1
\in A$ by $\circledast_0$.

By the definition of ${\Cal A}'$, for every $a \in{\Cal A}' \cap {\Cal
A}^*$, we have
\mr
\item "{$\circledast_6$}"  $n \in B_2 \Rightarrow 
{ \text{\rm Ch\/}}_a(n) \ge f^*(n) > f_{\delta(a)}(n)$
\sn
\item "{$\circledast_7$}"  $n \in B_1 \Rightarrow 
{ \text{\rm Ch\/}}_a(n) =
\cup\{f_{\beta_\varepsilon}(n):\varepsilon \in R_2 \cap \delta(a)\}
\le f_{\delta(a)}(n)$.
\ermn
But this contradicts the observation below.
\enddemo
\bigskip

\demo{\stag{nr.8} Observation}  If $B \subseteq \omega$, \ub{then} for
some $\varepsilon < \omega_1$ we have:

if $a \in {\Cal A}$ is $M^*$-closed and 
$\{n < \omega:\sup(a \cap \lambda_n) \le
f_{\sup(a)}(n)\} = B \text{ mod } J^{\text{bd}}_\omega$, \ub{then} otp$(a) <
\varepsilon$.
\enddemo
\bigskip

\demo{Proof}  Read the definition of ${\Cal A}$ 
(and ${\Cal A}^\varepsilon,{\Cal A}^\varepsilon_\delta$) and subclaim
\scite{nr.4} particularly $(*)_3$.
\nl
${{}}$  \hfill$\square_{\scite{nr.8}},
\square_{\scite{nr.7}},\square_{\scite{nr.1}}$
\enddemo
\bigskip

\remark{Remark}  Clearly 
\scite{nr.8} shows that we have much freeness in the choice of
${\Cal A}^\varepsilon_\delta$'s.
\endremark
\bn
We can get somewhat more, as in \cite{Sh:e}
\proclaim{\stag{nr.10} Claim}  In Claim \scite{nr.1} we can add to the
conclusion
\mr
\item "{$(*)$}"  ${\Cal A}$ satisfies the diamond,
i.e. $\diamondsuit_{\Cal A}$.
\endroster
\endproclaim
\bigskip

\demo{Proof}  In \scite{nr.4} we can add
\mr
\item "{$(*)_5$}"  $\{2n+1:n < \omega\} = \emptyset$ mod
$D_{\ell,\varepsilon}$ for $\ell < 2,\varepsilon < \omega_1$. 
\ermn
This is easy:  replace $D_{\ell,\varepsilon}$ by
$D'_{\ell,\varepsilon} = \{A \subseteq \omega:\{n:2n \in A\} \in
D_{\ell,\varepsilon}\}$.  We can fix a countable vocabulary $\tau$ and
for $\zeta < \omega_1$ choose a function $F_\zeta$ from 
${\Cal P}(\omega)$ onto $\{N:N$ is a $\tau$-model with universe $\zeta\}$ such
that $F_\zeta(A) = F_\zeta(B)$ if $A=B$ mod finite. 
\bn
\ub{Case 1}:  $\mu > 2^{\aleph_0}$.

Lastly, for $a \in {\Cal A}$ let $\delta_a,\varepsilon_a$ be such that $a
\in {\Cal A}^{\varepsilon_a}_{\delta_a}$ and let $A_a = \{n:\sup(a \cap
\lambda_{2n+1}) < f_{\delta_a}(2n)\}$, and let $N_a$ be the $\tau$-model with
universe $a$ such that the one-to-one order preserving function from
$\zeta$ onto $a$ is an isomorphism from $F_\zeta(N)$ onto $N$.
Note that in the proof of ``${\Cal A} \subseteq [\lambda]^{\aleph_0}$
is stationary", i.e. of \scite{nr.6}, given a $\tau$-model $M$ with
universe $\lambda$ \wilog \, $\lambda_0 > 2^{\aleph_0}$ and so can
demand that the isomorphism type of $M_\eta$ is the same for all $\eta
\in \text{ lim}(T)$ and, of course, $M \in M_\eta$.  Hence the
isomorphic type of $M \restriction u_\eta$ is the same for all $\eta
\in \text{ lim}(T)$ where $u_\eta$ is the universe of  $M_\eta$.  Now
in the choice for $B_1$ we can add the demand
$F_{\varepsilon^*}(\{n:2n+1 \in B_1\})$ is isomorphic to $M
\restriction u_\eta$ for every $\eta \in \text{ lim}(T)$.
\nl
Now check.   
\bn
\ub{Case 2}:  $\mu \le 2^{\aleph_0}$.

Similarly letting $\{2n +1:n < \omega\}$ be the disjoint union of
$\langle B^*_n:n <\omega\rangle$, each $B_n$ infinite.  We use
$A_a \cap B^*_n$ to code model with universe $\subseteq \zeta$ for
some $\zeta < \omega_1$, by a function $\bold F_n$.  We then
let $N_a$ be the model with universe $a$ sucht hat the order
preserving function from $a$ onto a countable ordinal $\zeta$ is an
isomorphism from $N_a$ onto $\cup\{\bold F_n(A_a \cap B^*_n):n <
\omega\}$ when the union is a $\tau$-model with universe $\zeta$.

Now we cannot demand them all $M_\eta,\eta \in \text{ lim}(T)$ has the
same isomorphism type but only the same order type.  The rest should
be clear.  \hfill$\square_{\scite{nr.10}}$
\enddemo
\bn
We can also generalize
\proclaim{\stag{nr.11} Claim}  We can weaken the assumption of
\scite{nr.1} as follows
\mr
\item "{$(a)$}"  $\lambda = \,\text{\rm cf}(\lambda) > \mu$ instead
$\lambda = \mu^+$ (still necessarily $\mu_* \le \mu$)
\sn
\item "{$(b)$}"   replace $J^{\text{bd}}_\omega$ by an ideal $J$ on $\omega$
containing the finite subsets, $\lambda_n = \,\text{\rm cf}(\lambda_n) >
\aleph_1,\mu = \,\text{\rm lim}_J \langle \lambda_n:n < \omega \rangle$ but
not necessarily $n < \omega \Rightarrow \lambda_n < \lambda_{n+1}$ and
add ${\Cal P}(\omega)/J$ is infinite (hence uncountable).
\endroster
\endproclaim
\bigskip

\demo{Proof}  In \scite{nr.4} in $(*)_3$ we choose $\langle
A_\varepsilon:\varepsilon < \omega_1 \rangle$, a sequence of subsets
of $\omega$ such that $\langle A_\varepsilon/J:\varepsilon < \omega_1
\rangle$ are pairwise distinct.  This implies some changes and waiving
$\lambda_n < \lambda_{n+1}$ requires some changes in \scite{nr.6}, 
in particular for each $n$ using $\langle \bold B_\alpha:\alpha \in
S^{\lambda_n}_{\aleph_0}\rangle$ with $\bold B_\delta = \{\eta \in
\text{ lim}(T^*):a \cap \lambda_n \subseteq \alpha\}$ and the
partition theorem \cite[XI,3.7,pg.549]{Sh:f}.
\nl
${{}}$    \hfill$\square_{\scite{nr.11}}$ 
\enddemo
\newpage


\nocite{ignore-this-bibtex-warning} 
\newpage
    
REFERENCES.  
\bibliographystyle{lit-plain}
\bibliography{lista,listb,listx,listf,liste}

\enddocument